\documentclass[11pt,a4paper]{article}
\usepackage{float}

\usepackage[margin=1in]{geometry}

\usepackage[compress]{cite}
\usepackage{color}
\usepackage{amsopn}
\usepackage{lipsum}
\usepackage{amsfonts}
\usepackage{graphicx}
\usepackage{epstopdf}
\usepackage{algorithmic}
\usepackage{multirow}  
\usepackage{verbatim}
\usepackage{appendix}
\usepackage{subcaption}
\usepackage{amsmath}
\ifpdf
  \DeclareGraphicsExtensions{.eps,.pdf,.png,.jpg}
\else
  \DeclareGraphicsExtensions{.eps}
\fi

\newenvironment{keywords}
{\par\medskip\noindent\textbf{Keywords: }\ignorespaces}
{\par\medskip}

\title{Large deviations, gradient flows and fluctuating hydrodynamics for fractional diffusion}

\author{Julien Barr\'e$^1$ and Hiroshi Horii$^{1,2}$}

\date{
\small
$^1$Institut Denis Poisson, Universit\'e d'Orl\'eans, CNRS,
Universit\'e de Tours\\
$^2$Universit\'e de la C\^ote d'Azur, CNRS,
Observatoire de la C\^ote d'Azur,
Boulevard de l'Observatoire,
CS 34229 06304 Nice Cedex 4, France\\[0.4em]
}

\begin{document}
\maketitle

\begin{abstract}
This work is a step towards a general understanding of fluctuations for systems of particles described at large scale by fractional diffusion. First, we derive formally a Large Deviation Principle for the empirical measure of mean-field interacting particles with a Lévy-type noise. In the case of independent particles, this unveils a generalized gradient flow structure at large scales. Second, we derive and solve numerically a stochastic Partial Differential Equation describing the fluctuating hydrodynamics of this system. The results compare well with direct numerical integration of the particles' dynamics. 
\end{abstract}


\section{Introduction}
\label{sec:intro}
Very different microscopic dynamics for a system of many particles may lead at large scales to qualitatively similar diffusive macroscopic dynamics for the particle density: Brownian particles, or particles undergoing Langevin dynamics, run and tumble, or other particles undergoing a jump process for velocities, particles ruled by an exclusion process, more general lattice gases... For a few references, \cite{solon2015active} provides a detailed comparison of active Brownian and run-and-tumble particles. \cite{MFT2015}, section 8, discusses many examples of lattice gas type. 
Crucially, this "universality" at large scales is not limited to the macroscopic deterministic equation, but encompasses the description of finite-size fluctuations: these systems often share a similar large deviation structure, and this underpins the Macroscopic Fluctuations Theory (MFT) \cite{MFT2002,MFT2015,derrida2007non}, a general theory of large-scale fluctuations for diffusive systems. This theory not only describes small fluctuations (of central limit theorem (CLT) type), but also large ones, at the large deviation level. 

The picture for macroscopic fluctuations for systems that are diffusive at large scales is hence very complete. In particular:
\begin{itemize}
\item The empirical measure of $N$ particles satisfies a dynamical large deviation principle, with a rate function $I(\rho)$ of the following type:
\begin{equation}
I[\rho] =\int_0^T \sup_\varphi \left[\int_\Omega \varphi \partial_t \rho dx -\mathcal{H}(\rho,\varphi)\right] dt,
\label{eq:Ratediff}
\end{equation}
where $\Omega$ is the spatial domain and the Hamiltonian
\begin{equation}
\mathcal{H}(\rho,\varphi)=\int_\Omega \big\{\varphi \big[\nabla(D(\rho)\nabla \rho) \big] +\chi(\rho)\nabla \varphi\cdot \nabla \varphi \big\}dx.
\label{eq:Hdiff}
\end{equation}
These expressions involve two density-dependent functions (which can be symmetric matrices in the multidimensional case): the diffusion $D(\rho)$, and the mobility $\chi(\rho)$.
 The rate function can be rewritten as
$I[\rho] = \frac14 \int_0^T||\partial_t \rho-\nabla(D(\rho)\nabla \rho)||^2_{-1,\chi(\rho)}{dt}$, where $||\cdot||_{-1,\sigma}$ is the $\sigma$-weighted negative Sobolev $H^{-1}$ norm.
\item All the microscopic information is summarized in these two functions, which reflects a high degree of universality with respect to the microscopic dynamics for the large-scale fluctuations of diffusive systems.
\item If the system is reversible, this large deviation principle is deeply related to a gradient flow structure \cite{JKO,adams2011large}; this gradient flow is defined by a free energy functional and a distance of Wasserstein$-2$ type. 
\item It is possible to write a Dean-Kawasaki (DK) equation\cite{Dean1996,Kawasaki1998}: this SPDE complements the macroscopic diffusive equation with a small multiplicative noise proportional to $N^{-\frac12}$; it is a natural candidate to simulate the fluctuating empirical measure at a mesoscopic level. It reads
\begin{equation}
\partial_t \rho = \nabla\cdot(D(\rho)\nabla \rho) +\nabla\cdot\left(\sqrt{\frac{2\chi(\rho)}{N}}\xi(x,t)\right),
\label{eq:DKdiff}
\end{equation}
where $\xi$ is a Gaussian random field, delta-correlated in space and time.

\item The DK equation \eqref{eq:DKdiff} is a singular nonlinear SPDE for which the notion of solution is very subtle, see, for instance, \cite{renesse,fehrman2019well,fehrman2024well}, which could raise questions on its relevance. However, discretizations of the Dean–Kawasaki equation, such as finite-difference or finite-element schemes, or regularizations of it, have been rigorously shown to approximate the density fluctuations of $N$ 
diffusing particles, either independent or interacting in a mean-field way
\cite{cornalba2023dean,djurdjevac2024weak,cornalba2026density}. These papers show that the DK equation is much more precise than the linear fluctuating hydrodynamics equation, which would correspond to the CLT, and an approximation of order $N^{-\frac12}$.  
DK equation hence provides a viable numerical tool to analyze mesoscopic fluctuations for diffusive systems.


\end{itemize}

Our goal in this article is to investigate the possible generalization of the points above for systems governed at large scales by fractional diffusion (i.e. with a fractional Laplacian instead of a standard one), including some tests of numerical computations. To gain some insight into this broad question, we shall study primarily a very simple system: $N$ particles undergoing independent Lévy flights, adding in some cases a mean-field interaction. 

We now review some related literature. Anomalous diffusion is a broad subject  (see, for instance, \cite{metzler2014anomalous} for a review), and the diffusion equation with a fractional Laplacian has found applications in various scientific domains. We may mention, among others: search strategies (see, for instance, \cite{raposo2009levy}), anomalous diffusion of cold atoms in optical lattices \cite{sagi2012observation,kessler2012theory}, or anomalous heat transport \cite{delfini2008nonequilibrium,lepri2009stochastic, priyanka2018anomalous,andreucci2025thermal} (see also, for instance, \cite{dhar2019anomalous} for a review). For this last example, rigorous proofs of super-diffusive behavior start from microscopic dynamics \cite{jara2015superdiffusion,bernardin20163, saito20195}. In this context, a model of independent Lévy walks has also been extensively studied \cite{cipriani2005anomalous,lepri2011density,dhar2013exact}; the solution in \cite{dhar2013exact} includes a computation of the large deviation function for the current. An alternative approach to anomalous heat conduction relies on nonlinear fluctuating hydrodynamics \cite{spohn2014nonlinear}, and leads to a description of fluctuations by a linear SPDE
\cite{kundu2019fractional}, bearing some similarities with our results (our Eq.\eqref{eq:TCL}). 
In addition, anomalous heat transport in systems with heavy-tailed microscopic dynamics has been studied in \cite{horii2022anomalous}, and the dynamical large-deviation function in related systems has been studied in \cite{horii2022large}.
Besides these applications, the question of generalizing MFT to systems governed at large scales by fractional diffusion, or a generalization thereof, has been explicitly attacked recently in \cite{Bernardin2025,bernardin_etal_2026}, which study long-range interacting lattice-based Ginzburg-Landau models. In this rather general setting, and considering non trivial boundary conditions driving the system genuinely out of equilibrium, the authors obtain at large scales a dynamical large deviation principle which is similar to the diffusive one, with inverse Sobolev norms generalizing the weighted $H^{-1}$ norms typical of diffusive systems. 
By contrast, we will show that the mathematical structure at large scales in our simple setting is different, featuring LDPs of exponential form characteristic of jump processes. \footnote{Long-range interacting Ginzburg-Landau models do present jump-like LDPs when the microscopic noise is Poissonian \cite{CedricRaphael}.}
Systems of mean-field interacting Brownian particles approximate McKean-Vlasov processes \cite{sznitman} and provide a very large class of models described by MFT-like large deviation principles; see \cite{Dawson1987} for an early derivation of such an LDP. MFT has been applied recently to the related Dyson gas, where singular long range interactions trigger anomalous diffusion\cite{dandekar2024current}.  It is natural to study in the same context generalized McKean-Vlasov processes, driven by Lévy instead of Brownian noise. While the law of large numbers has been studied for these processes \cite{jourdain2007nonlinear,bao2024note,cavallazzi2025well,tran2025infinite,cavallazzi2025quantitative}, and lead to some nonlinear fractional diffusion equation for the particle density, we are not aware of large deviation studies. We note that \cite{dandekar2023dynamical} derives from $N$ particles a fluctuating PDE involving a kind of nonlinear fractional diffusion, but the noise itself is the same as for Brownian particles.

This paper is organized as follows. Drawing on the literature on jump processes \cite{Erbar2014,Mielke2014,Peletier2022}, Section~2 establishes a central limit theorem and a dynamical large deviation principle (LDP) for the empirical measure of the $N$ particles Lévy–flight system and identifies the associated generalized gradient flow structure, which is not of Wasserstein type. We extend this LDP formally to the case of mean-field interacting particles. As a consequence of the LDP structure, the equivalent of the Dean-Kawasaki equation would have non-Gaussian noise, a priori very difficult to simulate. In Section~3, we then resort to i) the CLT (which can also be seen as a linearization of DK) or ii) a Gaussian truncation of DK (but keeping the multiplicative structure of the noise) in order to simulate the fluctuations of the system at mesoscale. 
We devise and implement a numerical method to solve the obtained SPDEs, based on the spatial Fourier transform. 
Direct comparisons show a very good agreement with microscopic simulations.

\section{Large deviation, gradient flow structure, and central limit theorem}
\label{sec:main}

For simplicity, we will stick to one space dimension and periodic boundary conditions. The spatial domain is hence $\Omega=[0,L]$, with $L=2\pi$. We consider first the simplest situation described macroscopically by a fractional diffusion equation: $N$ independent particles, each one undergoing a stable symmetric L\'evy process with parameter $\alpha \in ]0,2]$. We define the empirical measure
\[
\mu^N = \frac1N \sum_{i=1}^N \delta\big(x-X_i(t)\big), 
\]
where the $X_i$ are independent L\'evy processes such that, for all $s$
\begin{equation}
\label{eq:characteristic}
\mathbb{E}\left[e^{is \big(X_i(t)-X_i(0)\big)}\right] = e^{-K(\alpha) t |s|^\alpha},
\end{equation}
where \(K(\alpha)=\sigma^\alpha/2\) with \(\sigma\) the scale parameter. We recall that periodic boundary conditions are used for the $X_i(t)$. For \(\alpha=2\), this yields Gaussian increments with variance \(\sigma^2 t\). In particular, when \(\sigma=1\), the process \(X_i\) is a standard periodized Brownian motion. 
The initial positions $X_i(t=0)$ are drawn independently at random from the probability density $\rho_0$. By the law of large numbers, when $N$ tends to infinity $\mu^N$ tends to the measure with density $\rho(x,t)$, such that $\rho$ is solution to the following fractional diffusion equation
\begin{equation}
\partial_t \rho =\mathcal{L}_\alpha \rho, \quad \rho(x,t=0)=\rho_0(x),
\end{equation}
where $\mathcal{L}_\alpha$ is the fractional Laplacian with parameter $\alpha$. $\mathcal{L}_\alpha$ is also the infinitesimal generator associated to the L\'evy process.
We use Fourier series defined for a function $g$ as
\[
\hat{g}(k) = \frac{1}{2\pi} \int_{\Omega} g(x) e^{-ikx} dx;
\]
the operator $\mathcal{L}_\alpha$ then reads in Fourier space
\begin{equation}
\widehat{\mathcal{L}_\alpha \rho}(k) = -K(\alpha) |k|^\alpha \hat{\rho}(k).
\label{eq:fracL}
\end{equation}
This Fourier-series representation corresponds to the spectral definition of the fractional Laplacian in [\citenum{roncal2014transference}, Section 1]. In addition, see \cite{lischke2020fractional} for a thorough discussion of different possible definitions of the fractional Laplacian on a bounded domain.

\subsection{Large deviation principle for the empirical measure}
\label{sec:ldp}
$\mu^N(t)$ is a family of paths of random measures on $\Omega$. Using classical results (see for instance \cite{Feng1994} Theorem 3.1; for a formal derivation, see \cite{Feliachi2021} section 4.1), we know that $\mu^N$ satisfies a large deviation principle with speed $N$ and Hamiltonian
\begin{equation}
\mathcal{H}[\mu,\varphi] = \int_\Omega \mu(x) e^{-\varphi(x)}\mathcal{L}_\alpha[e^\varphi](x) dx.
\label{eq:Ham}
\end{equation}
The rate function $\mathcal{I}$ of the large deviation principle for the empirical density $\mu$ is then (we forget here about the initial condition contribution):
\begin{equation}
\mathcal{I}[\mu,0\leq t\leq T] = \int_0^T \left( \sup_\varphi \int_\Omega \varphi \partial_t \mu\, dx -\mathcal{H}[\mu,\varphi]\right) dt.
\label{eq:LDF}  
\end{equation}

\subsection{Large-deviation Hamiltonian for mean-field interacting Lévy particles}
\label{sec:mean-field-levy}
The previous large-deviation Hamiltonian was derived for independent symmetric Lévy flights. It is interesting to consider the analogous large-deviation Hamiltonian when the particles interact through a mean-field drift, which is the setting studied, for instance, in \cite{jourdain2007nonlinear,bao2024note,cavallazzi2025well,tran2025infinite,cavallazzi2025quantitative}. Take the system
\[
dX_j(t)=F_j^N(t)\,dt+dL_j^\alpha(t), \qquad F_j^N(t):=\frac1N\sum_{\ell\neq j} b(X_\ell(t)-X_j(t)),
\]
where \(b\) is a given function describing the interaction. For a short time interval \(\Delta t\), we write
\[
\Delta X_j=F_j^N(t)\Delta t+\Delta t^{1/\alpha}L_j^\alpha .
\]
The large deviation Hamiltonian is given formally by the computation:
\begin{equation}
\mathcal H_{\rm mf}[\mu,\varphi]:=\lim\limits_{\Delta t \to 0}\lim\limits_{N \to \infty}\frac{1}{N\Delta t}\log \mathbb E\left[e^{N\left\langle\mu^N(t+\Delta t)-\mu^N(t),\varphi\right\rangle}\right], \label{eq:Hdef}
\end{equation}
with
\[
N\left\langle\mu^N(t+\Delta t)-\mu^N(t),\varphi\right\rangle=\sum_{j=1}^N\left(\varphi(X_j(t+\Delta t))-\varphi(X_j(t))\right).
\]
The expectation in \eqref{eq:Hdef} is to be understood here as a conditional expectation at fixed $X_j(t)$, for all $j$.
Introducing the function $g(x)=e^{\varphi(x)}$, we have:
\begin{eqnarray*}
    \mathbb{E}\left[e^{\varphi(X_j+\Delta X_j)-\varphi(X_j)}\right] &=&
    g(X_j+\Delta X_j)e^{-\varphi(X_j)} \\
    &=& \mathbb{E}\left[ \sum_{k\in \mathbb{Z}} \hat{g}(k)e^{ikX_j}e^{ik \Delta X_j} \right] e^{-\varphi(X_j)} \\
    &=& \sum_{k\in \mathbb{Z}} \hat{g}(k)e^{ikX_j}e^{-K(\alpha)|k|^\alpha \Delta t} e^{i k F_j^N \Delta t}
    e^{-\varphi(X_j)} \\
    &=& 1+ \left(\mathcal{L}_\alpha[e^\varphi](X_j)e^{-\varphi(X_j)}\right) \Delta t +F_j^N \varphi^\prime(X_j) \Delta t +o(\Delta t)
\end{eqnarray*}
Inserting in \eqref{eq:Hdef}, and expanding the logarithm, we are left with
\begin{eqnarray}
\mathcal H_{\rm mf}[\mu,\varphi]&=&\lim\limits_{\Delta t \to 0}\lim\limits_{N \to \infty} \frac1N \sum_j \left(\mathcal{L}_\alpha[e^\varphi](X_j)e^{-\varphi(X_j)} \right. \nonumber \\
&& \left. +\frac1N \sum_{\ell\neq j}b(X_\ell-X_j) \varphi^\prime(X_j) +o(1)\right) \nonumber 
\end{eqnarray}
We finally obtain the large deviation Hamiltonian: 
\begin{equation}
\mathcal H_{\rm mf}[\mu,\varphi]
= \int_\Omega \mu(x) F_\mu(x) \partial_x \varphi(x) dx+\int_\Omega \mu(x) e^{-\varphi(x)} \mathcal{L}_\alpha[e^\varphi](x) dx,
\end{equation}
with
\[
F_\mu(x) =\int_{\Omega} b(y-x) \mu(y) dy.
\]
Hence, the interaction drift term simply adds a corresponding term in the Hamiltonian.
This computation is formal and calls for a full justification.

\subsection{Properties of the rate function}
The empirical measure process in \ref{sec:ldp} describes independent particles.  
Introducing the functional (minus the entropy)
\[
U[\mu] = \int_\Omega \mu(x) \ln \mu(x)  dx, 
\]
we expect that $U$ is the quasi-potential associated to the dynamical large deviation principle \eqref{eq:Ham}-\eqref{eq:LDF}. 
Indeed we can check that $\mathcal{H}$ satisfies $\mathcal{H}[\mu,-\varphi]= \mathcal{H}\left[\mu, \varphi+\frac{\delta U}{\delta \mu}\right]$, which is the detailed balance relation with quasi potential $U$:
\begin{eqnarray}
\mathcal{H}\left[\mu, \varphi+\frac{\delta U}{\delta \mu} \right] &=& \int_\Omega \mu \frac{e^{-\varphi}}{\mu}\mathcal{L}_\alpha\left[\mu e^{\varphi}\right] \\
&=&   \int_\Omega \mathcal{L}_\alpha \left[e^{-\varphi}\right] \mu e^{\varphi} \\
&=& \mathcal{H}[\mu,-\varphi],
\end{eqnarray}
where we have used that $\mathcal{L}_\alpha$ is self adjoint for the $L^2$ scalar product. In particular, the functional $U$ satisfies the stationary Hamilton-Jacobi equation
$\mathcal{H}\left[\mu,\frac{\delta U}{\delta \mu}\right]=0$.

In the interacting case of section \ref{sec:mean-field-levy}, it is not obvious what the quasi-potential is, nor whether it satisfies a detailed-balanced relation.

\subsection{Generalized gradient flow structure}
The reversible large deviation principle given by the Hamiltonian \eqref{eq:Ham} yields a generalized gradient flow structure, which, unlike classical gradient flow structures, requires no distance.
These structures have now been found for many macroscopic equations for which it can be related to an underlying microscopic system satisfying a large deviation principle; the microscopic dynamics is usually of jump type, see for instance \cite{Erbar2014,Mielke2014,Bouchet2020,Peletier2022}.
Like classical gradient flows, they can be put to use to prove hydrodynamic limit results \cite{Hoeksema2023}.

Introducing the dissipation function: 
\begin{eqnarray}
R^\ast[\mu,\varphi]&=& \mathcal{H}\left[\mu,\varphi+\frac12 \frac{\delta U}{\delta \mu}\right] - \mathcal{H}\left[\mu, \frac12 \frac{\delta U}{\delta \mu}\right]\\
&=& \int_\Omega \sqrt{\mu}e^{-\varphi}\mathcal{L}_\alpha[\sqrt{\mu}e^{\varphi}]dx -  \int_\Omega \sqrt{\mu} \mathcal{L}_\alpha[\sqrt{\mu}]dx,
\end{eqnarray}
the fractional diffusion equation can then be written as 
\begin{equation}
\partial_t \mu =  \frac{\delta R^\ast}{\delta \varphi} \left[\mu,-\frac12 \frac{\delta U}{\delta \mu} \right].
\end{equation}
At this point, we may comment further on \cite{Bernardin2025}: the authors start from a microscopic model whose deterministic hydrodynamic limit is a nonlinear fractional diffusion equation. They derive a large deviation principle, and it is worthwhile to notice that they find a different structure than the jump-like one we have here: their large deviation rate function is quadratic (see their section 4.1), and their deterministic equation then appears as a "classical" gradient flow.

\subsection{A stochastic PDE describing the fluctuations?}
As explained in the introduction, for standard diffusive systems, the finite-size fluctuations of the empirical measure are described by the Dean-Kawasaki equation \eqref{eq:DKdiff}, a nonlinear stochastic PDE (SPDE) with  Gaussian noise. This SPDE can 
in some cases be formally derived directly from the particle dynamics \cite{Dean1996}. It can also be guessed from the large deviation principle \eqref{eq:Ratediff}-\eqref{eq:Hdiff}; the Gaussian character of the noise in \eqref{eq:DKdiff} is related to the fact that the large deviation Hamiltonian \eqref{eq:Hdiff} is quadratic in the momentum variable.

For fractional diffusion, we have seen that the Hamiltonian \eqref{eq:Ham} is not quadratic in the momentum $\varphi$, and has the structure of the Hamiltonian of a jump process rather than of a diffusion. It is still in principle possible to write down an SPDE corresponding to the large deviation principle, but the noise would not be Gaussian, and the equation would be very difficult to handle analytically or numerically. For these reasons, we have not explored this direction further. Alternatively, we can expand the Hamiltonian \eqref{eq:Ham} up to quadratic order in 
the momentum $\varphi$; this amounts to considering the fluctuations at the level of the Central Limit Theorem (CLT) only. 

This expansion of \eqref{eq:Ham} up to quadratic order in $\varphi$ is:
\begin{equation}
\mathcal{H}_2[\mu,\varphi] = \int_\Omega \mu(x) \mathcal{L}_\alpha[\varphi](x) dx +\frac12 \int_\Omega \mu(x) \Gamma_2^\alpha[\varphi](x) dx + O(\varphi^3),
\label{eq:Ham2}
\end{equation}
where $\Gamma_2^\alpha$ is the carr\'e du champ operator for the Levy process:
\begin{equation}
\Gamma_2^\alpha[\varphi](x) = \mathcal{L}_\alpha[\varphi^2](x) - 2 \varphi(x) \mathcal{L}_\alpha[\varphi](x).
\label{eq:cdc}
\end{equation}

In the next paragraph, we shall directly derive the CLT for our system of non-interacting Lévy flights, in the form of a linear stochastic PDE describing the statistics of the fluctuation of the empirical measure, and make the connection with the previous expansion \eqref{eq:Ham2}.

\subsection{Central limit theorem and linear stochastic PDE}
\label{sec:CLT}
Starting from a configuration at time $t$ given by the positions $(X_i)_{1\leq i\leq N}$, we consider a small time interval $\Delta t$. The $X_i(t)$ are independent stable symmetric Lévy processes, hence, for any $i$:
\begin{equation}
X_i(t+\Delta t) =X_i(t) +\Delta X_i~,~\Delta X_i = (\Delta t)^{\frac{1}{\alpha}} L_{i,t}^\alpha,
\label{eq:levy_flight}
\end{equation}
where the $L_{i,t}^\alpha$ are independent and drawn from the symmetric $\alpha-$ stable law, such that 
$\mathbb{E}[e^{isL_{i,t}^\alpha}]=e^{-K(\alpha)|s|^\alpha}$. We also consider the empirical measure
\[
\mu^N(t) = \frac1N \sum_{i=1}^N \delta(x-X_i(t))
\]
and are interested in its fluctuations. We take a test function $\varphi$ and compute the quantity:
\begin{eqnarray}
 \langle \mu^N(t+\Delta t),\varphi\rangle-\langle \mu^N(t),\varphi\rangle &=& \frac1N \sum \Delta \varphi_i,
 \end{eqnarray} 
where $\Delta \varphi_i = \varphi(X_i(t+\Delta t)) -\varphi(X_i(t))$.  To write a CLT, we need to compute
$\mathbb{E}[\Delta \varphi_i|X_i(t)]$ and $\mathbb{E}[(\Delta \varphi_i)^2|X_i(t)]$, conditional expectations given $X_i(t)$. 
Using the representation in Fourier space:
\begin{equation}
\Delta \varphi_i = \sum_{k\in \mathbb{Z}} \hat{\varphi}(k)e^{ikX_i} \left(e^{ik\Delta X_i}-1\right),
\end{equation} 
we have
\begin{equation}
\mathbb{E}[\Delta \varphi_i|X_i(t)] = \sum_{k\in \mathbb{Z}} \hat{\varphi}(k)e^{ik X_i}\left(e^{-K(\alpha)|k|^\alpha \Delta t}-1\right),
\end{equation}
where we have used the characteristic function of stable laws: $\mathbb{E}[e^{ikL^\alpha}] =e^{-K(\alpha)|k|^\alpha}$, and the scaling factor $K(\alpha)$ has been introduced above.
Keeping only terms of order $\Delta t$, and summing over the Fourier modes $k$:
\[
\mathbb{E}[\Delta \varphi_i] =\Delta t ~\mathcal{L}^\alpha[\varphi] (X_i).
\]
Hence
\begin{equation}
\frac1N \sum_i \Delta \varphi_i =   \Delta t ~ \left\langle \mu^N, \mathcal{L}^\alpha[\varphi] \right\rangle =  \Delta t ~ \left\langle \mathcal{L}^\alpha[\mu^N],\varphi \right\rangle.
\end{equation}
Similarly, we have
\begin{eqnarray}
\mathbb{E}[(\Delta \varphi_i)^2] &=& \sum_{k_1,k_2\in \mathbb{Z}} \hat{\varphi}(k_1)\hat{\varphi}(k_2)e^{i(k_1+k_2) X_i}\left(e^{-K(\alpha)|k_1+k_2|^\alpha \Delta t} \right. \nonumber \\
&& \left.-e^{-K(\alpha)|k_1|^\alpha \Delta t}-e^{-K(\alpha)|k_2|^\alpha \Delta t}+1\right).
\end{eqnarray}
Expanding in $\Delta t$, and summing over $k_1,k_2$:
\[
\mathbb{E}[(\Delta \varphi_i)^2]   = \Delta t ~\Gamma^\alpha[\varphi](X_i) +O(\Delta t^2),
\]
where the carr\'e du champ operator $\Gamma_2^\alpha[\varphi](x)$ is defined in \eqref{eq:cdc}.
Clearly $(\mathbb{E}[\Delta \varphi_i])^2 = O(\Delta t^2)$ can be neglected in the computation of the variance of $\Delta \varphi_i$.
Putting everything together, one obtains a CLT for the variation of the empirical measure between $t$ and $t+\Delta t$:
\begin{equation}
\langle \Delta \mu^N, \varphi\rangle = \langle \mathcal{L}^\alpha\mu^N,\varphi\rangle \Delta t + \frac{\sqrt{\langle \mu^N,\Gamma_2^\alpha[\varphi]\rangle \Delta t}}{\sqrt{N}} Z +O\left(\frac1N\right),
\label{eq:discrDK}
\end{equation}
where $Z$ is a standard normal random variable.

The expectation $\mu$ of the empirical measure is the expected solution of $\partial_t \mu =\mathcal{L}^\alpha \mu$. From \eqref{eq:discrDK}, one can write an equation for the fluctuations of the empirical measure $\eta =\sqrt{N}(\mu^N-\mu)$: for any test function $\varphi$:
\begin{equation}
\langle \Delta \eta,\varphi\rangle = \Delta t ~\langle\mathcal{L}^\alpha \eta, \varphi \rangle + \sqrt{\Delta t \langle \mu, \Gamma_2^\alpha[\varphi] }\rangle Z.
\end{equation}
This equation can be rewritten as a field equation, without a test function
\begin{equation}
\Delta \eta(x,t) =  \Delta t ~ \mathcal{L}^\alpha \eta (x,t) +  \sqrt{\Delta t}~ \xi(x,t),
\label{eq:discr_fh}
\end{equation}
where $\xi(x,t)$ is a random centered Gaussian field such that 
\[
\mathbb{E}[\xi(x,t)\xi(x',t')] = \delta(t-t') C(x,x').
\]
The correlation $C$ is determined by the identity, valid for any test function $\varphi$:
\[
\iint \varphi(x)\varphi(x') C(x,x') dxdx' = \int \mu(x) \Gamma_2^\alpha[\varphi](x) dx.
\]
In addition, in Fourier space, the correlation $C$ is described by using the Fourier representation of the carré du champ operator, which gives 
\begin{equation}
    \int \mu(x)\Gamma_2^{\alpha}[\varphi](x)=2\pi K(\alpha)\sum_{k_1,k_2\in\mathbb{Z}} \hat{\mu}(-k_1-k_2)\left[-|k_{1}+k_{2}|^{\alpha}+|k_{1}|^{\alpha}+|k_{2}|^{\alpha}\right]\hat{\varphi}(k_1)\hat{\varphi}(k_2).
\end{equation}
From this identity, we obtain the Fourier-space covariance $\hat{C}$:
\begin{equation}
\hat{C}(k_1,k_2)=\frac{K(\alpha)}{2\pi}\hat{\mu}(k_1+k_2)\left[-|k_{1}+k_{2}|^{\alpha}+|k_{1}|^{\alpha}+|k_{2}|^{\alpha}\right].
\label{eq:noise_fourier}
\end{equation}
The equation for the fluctuation $\eta$ can then be rewritten as a linear SPDE: 
\begin{equation}
\partial_t \eta = \mathcal{L}^\alpha \eta + \xi(x,t).
\label{eq:TCL}
\end{equation}
The noise is delta-correlated in time, but has non-trivial spatial correlation, which depends, of course, on the density $\mu$. 
\eqref{eq:TCL} is an important result because, coupled with the fractional diffusion equation for $\mu$ and after regularization of the noise by standard means, it gives a practical way to simulate the finite-size effects on the empirical density of particles. 
We note \cite{kundu2019fractional} obtains a similar linear fluctuating hydrodynamical equation as Eqs.\eqref{eq:noise_fourier}-\eqref{eq:TCL}, for $\alpha=3/2$ and a uniform temperature (ie a uniform density $\mu$ in the present setting): see their Eq. (32). When specializing \eqref{eq:noise_fourier} to $\alpha=3/2$ and a spatially uniform $\mu$, it recovers the same spectral structure as in \cite{kundu2019fractional}.

We also note that an alternative option to \eqref{eq:TCL} could be to solve the dynamical Lyapunov equation describing the evolution of the two-point correlation function $S(x,y,t)=\mathbb{E}\left[\eta(x,t)\eta(y,t)\right]$. A straightforward computation using Ito formula yields the equation:
\begin{equation}
\partial_t S(x,y,t) = \mathcal{L}^{\alpha,1}S(x,y,t)+  \mathcal{L}^{\alpha,2}S(x,y,t) +C(x,y),
\end{equation}
where $\mathcal{L}^{\alpha,1}$ (resp. $\mathcal{L}^{\alpha,2}$) acts on the first (resp. second) variable of the following two-variable function.
We have not followed this route further.

It is tempting to write from \eqref{eq:discrDK} a Dean-Kawasaki-like equation, with multiplicative noise, for the empirical density $\mu^N$, which would read
\begin{equation}
\partial_t \mu^N = \mathcal{L}^\alpha \mu^N + \frac{1}{\sqrt{N}}\xi[\mu^N](x,t),
\label{eq:DK_test}
\end{equation}
where $\xi[\mu^N](x,t)$, conditioned on $\mu^N(,\cdot,t)$, would be a Gaussian field, delta-correlated in time and with a non-trivial spatial correlation, depending on $\mu^N(\cdot,t)$ \footnote{Different interpretations of the multiplicative noise (Ito or Stratonovich) in \eqref{eq:DK_test} should differ only by an $O(1/N)$ term.}. 
For standard diffusion, the Dean-Kawasaki equation provides a better approximation to the true empirical density than the CLT \cite{Cornalba2019,cornalba2026density}.
In the present setting, it is not clear, however, that \eqref{eq:DK_test} would provide a better approximation of the true empirical density than \eqref{eq:TCL}, which has additive noise.

\section{Numerical simulations}
Equation~\eqref{eq:TCL} represents the linear fluctuating hydrodynamics for the system of L\'evy particles. Together with the fractional diffusion equation for the mean density, it provides a mesoscopic description of the fluctuating particle density at the level of the CLT. In this section, we present simulations of the stochastic PDE models \eqref{eq:TCL} and \eqref{eq:DK_test}, and compare them with direct simulations of a large number of particles. We do not aim at computational efficiency here; rather, our goal is to show that mesoscopic simulations can accurately reproduce the fluctuations of L\'evy flight particles undergoing fractional diffusion.
\subsection{Noise structure}
To simulate the stochastic PDE \eqref{eq:DK_test} we employ Gaussian noise whose correlation is specified in Fourier space \eqref{eq:noise_fourier}. Since the Fourier coefficients of a real–valued random field are generally complex, we decompose the noise as
\[
\hat{\xi}=\hat{\xi}_{\mathrm{Re}}+i\,\hat{\xi}_{\mathrm{Im}}.
\]
More generally, a complex $n$-component Gaussian vector $\mathbf{z}\in\mathbb{C}^{n}$ can be written as
$\mathbf{z}=\mathbf{x}+i\mathbf{y}$, where $\mathbf{x},\mathbf{y}\in\mathbb{R}^{n}$. 
In our setting, we denote the real and imaginary parts of the complex
Fourier space noise by
\[
\mathbf{x}\equiv\hat{\xi}_{\mathrm{Re}},
\qquad
\mathbf{y}\equiv\hat{\xi}_{\mathrm{Im}}.
\]
Then the complex random variable is represented in terms of its real and imaginary parts, with the following covariance matrix:
\[
\begin{bmatrix}
	\textbf{x} \\
	\textbf{y}
\end{bmatrix}
=\begin{bmatrix}\hat{\xi}_{\mathrm{Re}}\\\hat{\xi}_{\mathrm{Im}}\end{bmatrix}
\sim \mathcal{N} \left(
\boldsymbol{\mu}_{2n \times 1} =
\begin{bmatrix}
	\boldsymbol{\mu}_{\mathbf{x}} \\
	\boldsymbol{\mu}_{\mathbf{y}}
\end{bmatrix},
\boldsymbol{\Sigma}_{2n \times 2n} =
\begin{bmatrix}
	\boldsymbol{\Sigma}_{\mathbf{xx}} & \boldsymbol{\Sigma}_{\mathbf{xy}} \\
	\boldsymbol{\Sigma}_{\mathbf{yx}} & \boldsymbol{\Sigma}_{\mathbf{yy}}
\end{bmatrix}
\right),
\]
where $\boldsymbol{\Sigma}_{2n \times 2n}$ is a covariance matrix for this random variable. Hence, the complex Fourier‑space noise is sampled as a multivariate complex normal random variable.
In order to obtain the covariance matrix expression, we use the noise correlation obtained in Section \ref{sec:CLT}. The details of the computation are described in the appendix. By using this covariance matrix, we can generate a multivariate Gaussian random variable and simulate the stochastic PDEs \eqref{eq:TCL} or \eqref{eq:DK_test}.

\subsection{One-step validation of the SPDE model}
\label{sec:numerical_PDE}
In this section, we simulate the stochastic PDEs and compare them with the particle simulation and the theoretical result, which can be computed exactly in some cases. For this purpose, we use different test functions and different parameters $\alpha$. Recalling \eqref{eq:DK_test}, the approximated fractional stochastic diffusion equation is given by
\[
	\partial_t \mu (x,t) = \mathcal{L}^\alpha \mu(x,t) + \frac{1}{\sqrt{N}}\xi[\mu](x,t),
\]
where $\mu$ is the density profile. For \eqref{eq:TCL}, one uses the deterministic density profile instead of the fluctuating one. In the numerical simulations, we follow the same convention as in the theoretical part, which is adopted throughout the paper, namely \(K(\alpha)=\sigma^\alpha/2\) for every $\alpha$, and set $\sigma=1$. In addition, using an inner product with a test function $\varphi$, we can rewrite this equation, with 
$\xi^{\varphi}[\mu]:= \langle \xi[\mu],\varphi\rangle$):
\begin{equation}
	\langle \partial_t\mu(x,t) ,\varphi(x) \rangle =\langle \mathcal{L}^{\alpha} \mu(x,t) ,\varphi(x) \rangle +\frac{1}{\sqrt{N}}\xi^{\varphi}[\mu] (t)
\end{equation}
We use as initial condition for \eqref{eq:TCL} or \eqref{eq:DK_test} the periodized Gaussian
\begin{equation}
	\mu(x,0) = A_1 \left(e^{-(x-L/2)^2} + 1\right),
    \label{eq:IC}
\end{equation}
	with \(L = 2\pi\) and $A_1$ is a normalization factor. 
To validate the SPDEs, we adopt the following simulation process:
\begin{enumerate}
	\item \textbf{Single Time-Step Calculation:} \\
	First, perform a numerical calculation for a single time step $\Delta t$ and compute the variance of the inner product of $\Delta \mu$ and a test function $\varphi$ at $t=\Delta t$; this variance corresponds to the noise variance:
\[
\frac{\Delta t}{N}\,\mathbb{E}\Bigl[ (\xi^{\varphi}[\mu])^2 \Bigr]=\mathbb{E}\Bigl[\langle\Delta\mu(\Delta t),\varphi\rangle^2\Bigr]
-\Bigl(\mathbb{E}\bigl[\langle\Delta\mu(\Delta t),\varphi\rangle\bigr]\Bigr)^2.
\]
	
	\item \textbf{Monte Carlo Simulation:} \\
	Perform many independent one-step Monte Carlo simulations and average the resulting variances.
	
	\item \textbf{Comparison with Theoretical computation:} \\
	Compare the numerical experimental results with the theoretical values for each test function. 
	In this simulation, we use 
	\[
	\varphi(x) = \sin x \quad \text{and} \quad \varphi(x) = \cos x.
	\]
\end{enumerate}
We next perform the same one-step test directly using the microscopic
Lévy-flight dynamics. For this purpose, we shall compute the same quantity, namely  $\operatorname{Var}\bigl(\,\langle\Delta\mu_N,\varphi\rangle\bigr)=\frac{\Delta t}{N}\mathbb{E}[(\xi^{\varphi})^{2}]$ directly from the microscopic dynamics. Recalling the stochastic process \eqref{eq:levy_flight},
\[
	X_i(t+\Delta t) = X_i(t)+(\Delta t)^{1/\alpha}L_{\alpha},
\]
the variation $\Delta \mu_N$ is computed as:
\begin{equation}
	\langle \Delta\mu_N,\varphi \rangle = \frac{1}{N}\sum^N_{i=1}(\varphi(X_i(t+\Delta t))-\varphi(X_i(t))).
\end{equation}
When $t=0$, this expression is given by
\begin{equation}
	\langle \Delta\mu_N,\varphi \rangle = \frac{1}{N}\sum^N_{i=1}(\varphi(X_i(\Delta t))-\varphi(X_i(0))),
\end{equation}
To check the variances in the Lévy flight simulation, we adopt the following simulation process:
\begin{enumerate}
  \item \textbf{Initial PDF Setup and Particle Sampling:} \\
We determine the initial positions of the particles using the inverse Cumulative Distribution Function (CDF) method, with the probability distribution \eqref{eq:IC}.
\item \textbf{Simulating Particle Motion:} \\
For a small time step \(\Delta t\), each particle's position is updated according to a Lévy flight, resulting in a new position \(X(\Delta t)\). We then compute the change $\Delta \varphi$ from the initial position \(X(0)\) using a test function \(\varphi(x)\) by computing
\[
\Delta \varphi_i = \varphi\bigl(X_i(\Delta t)\bigr) - \varphi\bigl(X_i(0)\bigr).
\]
Also, we calculate the following
\[
\langle \Delta \mu_N,\varphi(x) \rangle=\frac{1}{N}\sum_i \Delta \varphi_i.
\]

\item \textbf{Variance Calculation:} \\
By repeating the particle update over many independent trials using a Monte Carlo simulation approach, the simulation computes statistical measures such as the mean and variance of $\langle \Delta \mu_N,\varphi(x) \rangle$. Note that the initial particles' positions $X_i(0)$ are fixed.
\end{enumerate}
In the $\alpha$‑stable law, the size of the jump is set by the scale parameter, here denoted $K(\alpha)$. Each Lévy increment is therefore multiplied by the factor $K(\alpha)^{1/\alpha}$. In addition, in agreement with the analysis in \cite{cornalba2023dean}, we choose the particle number $N$ and the spatial grid size $h$ such that the average number of particles per cell is sufficiently larger than one: $h\gg N^{-1}$.

Finally, to validate the above simulations, we compute the theoretical expression of the noise field correlation, which is given by
\begin{equation}
	\mathbb{E}[\xi^{\varphi}\xi^{\varphi}]=\mathbb{E}\left[\left(\int dx \xi(x) \varphi(x)\right)^2\right].
\end{equation}
To compute this correlation function, we use the Fourier series representation, which is given by
\[
	\xi^{\varphi}=\int \xi(x) \sum_{k\in\mathbb{Z}}\hat{\varphi}(k){e^{ikx}}dx.
\]
By substituting the above expression into the definition of the correlation function, when $\varphi(x)=\cos x$, the correlation function is
\begin{eqnarray}
	\mathbb{E}[\xi^{\varphi}\xi^{\varphi}]|_{\varphi=\cos x}=\pi^2\mathbb{E}[\hat{\xi}^2(k=1)+\hat{\xi}^2(k=-1)+2\hat{\xi}(k=1)\hat{\xi}(k=-1) ].
\end{eqnarray}
Recalling the Fourier space covariance of the noise given in \eqref{eq:noise_fourier}, we have
\[
	\mathbb{E}[\hat{\xi}(k_1)\hat{\xi}(k_2)]=\frac{1}{4\pi}\hat{\mu}(k_1+k_2)\left(-|k_1+k_2|^\alpha +|k_1|^\alpha+|k_2|^\alpha\right).
\]
Finally, we can obtain the following expression for the case of $\varphi(x)=\cos x$.
\begin{eqnarray}
	\mathbb{E}[\xi^{\varphi}\xi^{\varphi}]|_{\varphi=\cos x}=\pi\left(\hat{\mu}(0)+\frac{-2^{\alpha}+2}{4} \left(\hat{\mu}(2)+\hat{\mu}(-2) \right)  \right)
\end{eqnarray}
Note that the Fourier coefficients are defined by
\[
\widehat\phi(k)=\frac{1}{2\pi}\int_0^{2\pi}\phi(x)e^{-ikx}\,dx, \qquad k\in\mathbb Z.
\]
Therefore, for a normalized density, we have
\[
\hat{\mu}(0)=\frac{1}{2\pi}.
\]
As with this computation, we can also derive the correlation function when $\varphi(x)=\sin x$. It is given by
\begin{eqnarray}
	\mathbb{E}[\xi^{\varphi}\xi^{\varphi}]|_{\varphi=\sin x}=\pi\left(\hat{\mu}(0)-\frac{-2^{\alpha}+2}{4} \left(\hat{\mu}(2)+\hat{\mu}(-2) \right)  \right).
\end{eqnarray}
By substituting the numerically evaluated values $\hat{\mu}(0)$, $\hat{\mu}(2)$, and $\hat{\mu}(-2)$ into the above expressions, we obtain the theoretical values of the correlation function shown in Table~\ref{tab:combined}.
Table \ref{tab:combined} compares the analytical computation, the SPDE, and the direct particle simulation. The good agreement provides a sanity check: our simulation framework's random variable generation step is appropriate for the stochastic PDEs \eqref{eq:TCL}-\eqref{eq:DK_test}.

\begin{table}[ht]
  \centering
  \begin{tabular}{|c|c|c|c|c|}
    \hline
    $\varphi(x)$ & $\alpha$ & Theoretical & SPDE simulation  & Particle simulation  \\ \hline
    \multirow{4}{*}{$\cos x$}
      & $2.0$ & 0.460 & 0.457 & 0.451 \\ \cline{2-5}
      & $1.5$ & 0.483 & 0.488 & 0.487 \\ \cline{2-5}
      & $1.0$ & 0.500 & 0.505 & 0.497 \\ \cline{2-5}
      & $0.5$ & 0.512 & 0.523 & 0.511 \\ \hline
    \multirow{4}{*}{$\sin x$}
      & $2.0$ & 0.540 & 0.538 & 0.540 \\ \cline{2-5}
      & $1.5$ & 0.517 & 0.515 & 0.515 \\ \cline{2-5}
      & $1.0$ & 0.500 & 0.507 & 0.496 \\ \cline{2-5}
      & $0.5$ & 0.488 & 0.495 & 0.476 \\ \hline
  \end{tabular}
  \vspace{0.1cm}
  \caption{%
Variance of the test quantity $\langle\Delta\mu_N,\varphi\rangle$ for one time step at $\alpha=0.5,1.0,1.5,2.0$. In this table, (i) the theoretical prediction $E[(\xi^{\varphi})^{2}]$, (ii) the numerical result obtained with the noise \eqref{eq:noise_fourier} (1000 samples),
and (iii) an independent Monte-Carlo simulation of Lévy-flight dynamics (10\textsuperscript{5} samples). All entries are reported without the overall factor $\Delta t/N$. Since $N=10^3$ and $\Delta t=10^{-4}$ we have $\Delta t/N = 10^{-7}$, so the variances are $10^{-7}$ times the numbers listed. For the Lévy dynamics, we adopt the scale convention $K(\alpha)=\sigma^{\alpha}/2$ with $\sigma=1$, thus each Lévy increment is multiplied by $1/2^{1/\alpha}$. All three computations coincide, within error bars.}
  \label{tab:combined}
\end{table} 
For this one-step simulation, the additive noise and multiplicative noise models have the same noise covariance, since both are evaluated at the same initial density. The one-step noise statistics of the SPDE models have now been validated, and we next compare them with the microscopic particle dynamics over longer time scales.

\subsection{Long time behavior of fractional SPDE in numerical simulation}
\label{sec:longtime}
In the previous section, we confirmed that the SPDEs \eqref{eq:TCL}-\eqref{eq:DK_test} reproduce the short-time statistics of Lévy-flight particle systems and is consistent with the theoretical computation. In this section, we examine the long-time behavior of the SPDEs by comparison with particle simulations. 
Equation \eqref{eq:TCL} includes noise whose covariance depends on the deterministic density. In contrast, equation \eqref{eq:DK_test} incorporates multiplicative noise that directly depends on the stochastic density field itself. We shall also evaluate the differences between these two models.
\captionsetup[subfigure]{labelformat=parens,labelsep=space}
\begin{figure}[tbp]
  \centering
  \begin{subfigure}[t]{0.48\textwidth}
    \centering
    \includegraphics[width=\linewidth]{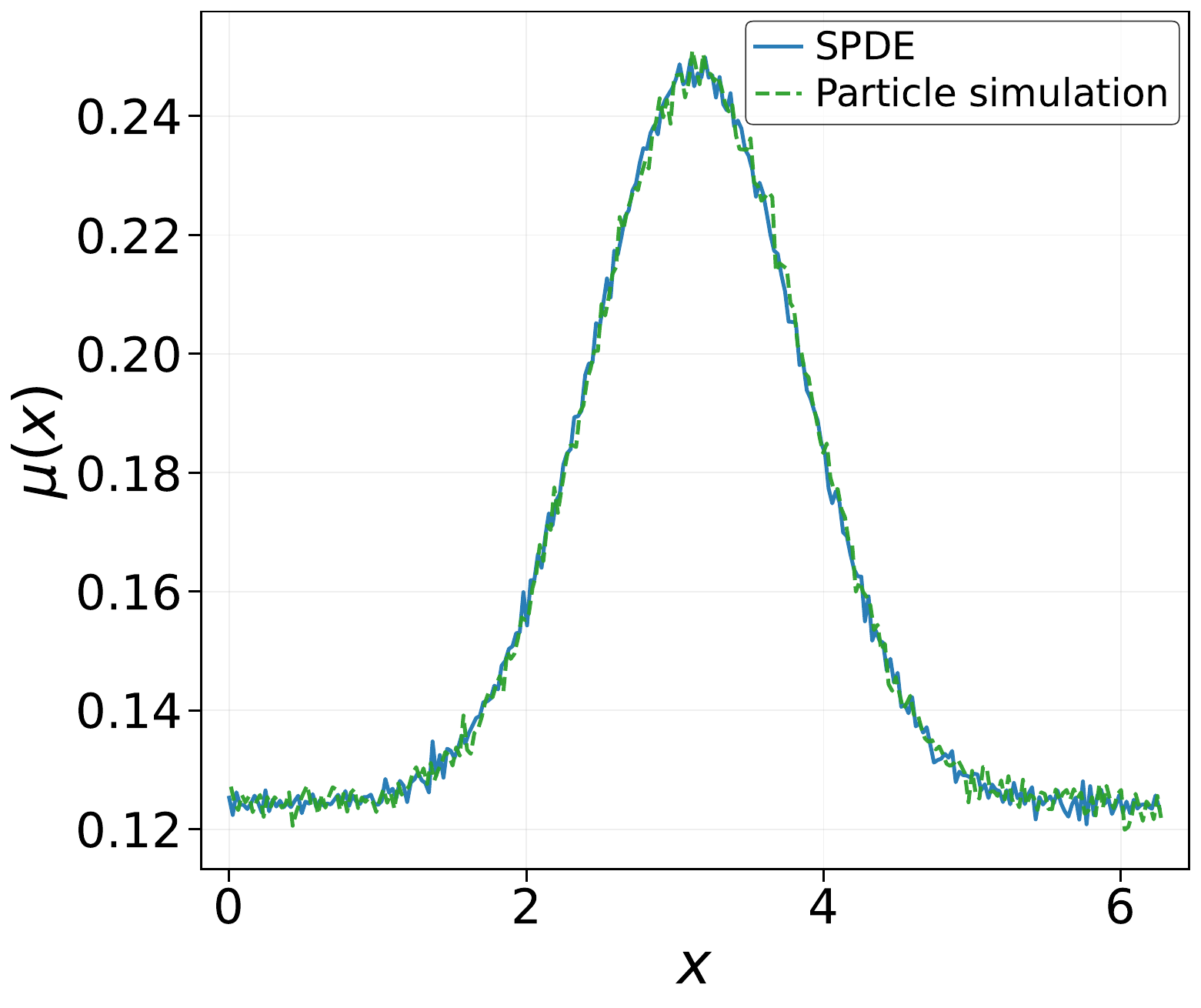}
    \subcaption{}
    \label{fig:density_profile_step1_random_initial}
  \end{subfigure}\hfill
  \begin{subfigure}[t]{0.48\textwidth}
    \centering
    \includegraphics[width=\linewidth]{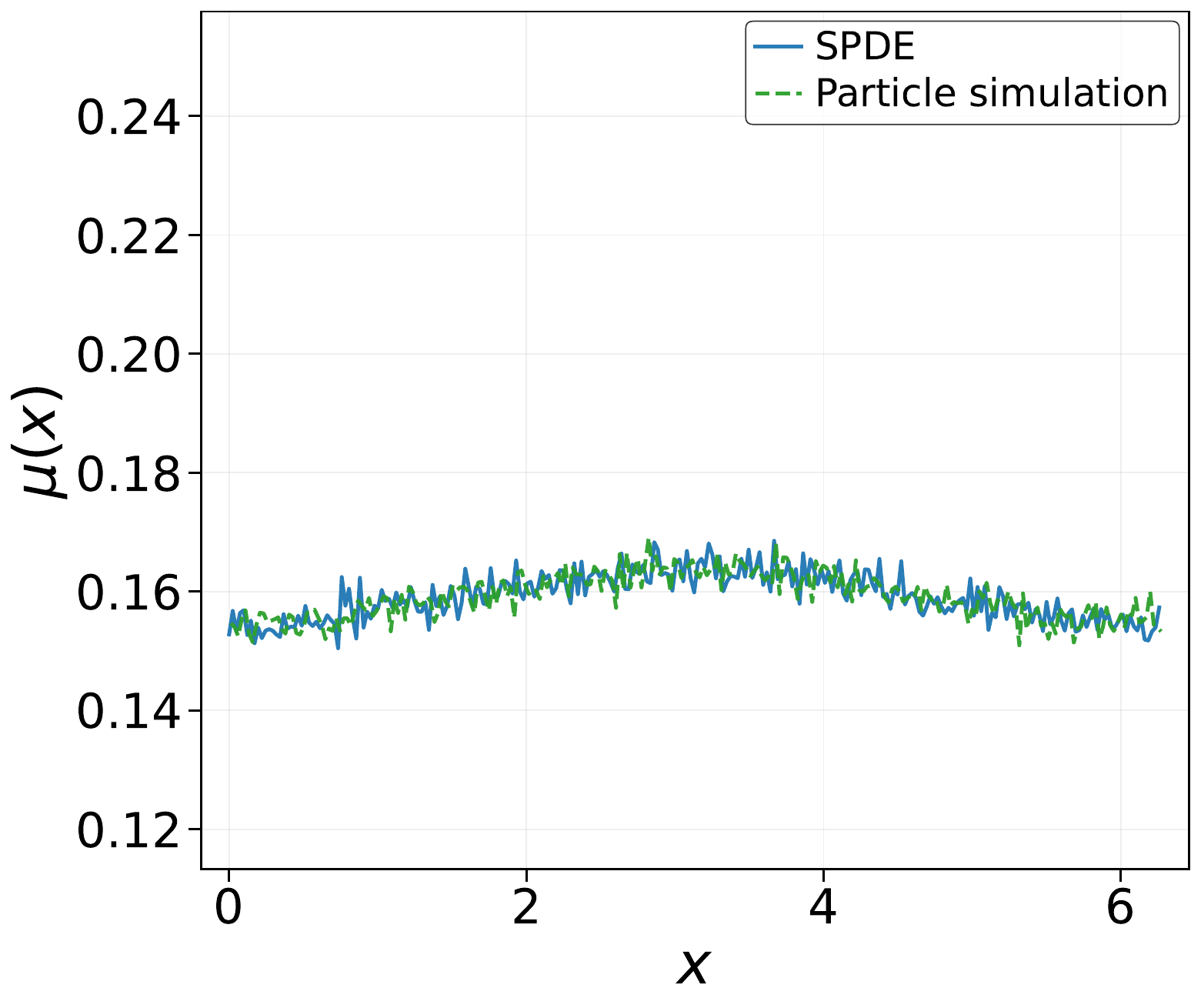}
    \subcaption{}
    \label{fig:density_profile_step500_random_initial}
  \end{subfigure}

  \caption{Representative density profiles for the fractional SPDE with additive noise and the particle simulation, for $\alpha=1.0$, $dt=0.01$, $N=2\times10^6$, and a spatial grid with 257 points. Panel (a) shows the profile at the first time step $t=0.01$, and panel (b) at $t=5.0$. The blue solid line denotes the SPDE with additive noise, and the green dashed line denotes the particle simulation. These plots are the basic time evolution of the profile. A histogram with 257 bins represents the particle density.}
  \label{fig:density_profiles_representative_random_initial}
\end{figure}
Before studying the fluctuation statistics, we first illustrate the time evolution of the density profile itself. Fig.~\ref{fig:density_profiles_representative_random_initial} shows representative profiles for the additive-noise SPDE, and the particle simulation in the case $\alpha=1.0$, $N=2\times10^6$, with 257 spatial grid points for the SPDE, which is solved by FFT. At the final time $t=5.0$, the profile has spread and become flatter. Following \cite{cornalba2023dean}, for quantitative comparisons we do not evaluate point-wise statistics of the density field, but rather its projection onto a fixed smooth test function \(\varphi\). Let \(\langle \,\cdot\,,\,\cdot\,\rangle\) denote the inner product. We define
\[
  C_1(t) \;=\; \mathbb{E}\!\big[\langle \mu(t), \varphi\rangle\big],\qquad
  C_2(t) \;=\; \mathbb{E}\!\Big[\big(\langle \mu(t), \varphi\rangle - C_1(t)\big)^2\Big],
\]
where $C_1$ is the mean value and $C_2$ is the variance. 
\captionsetup[subfigure]{labelformat=parens,labelsep=space}
\begin{figure}[tbp]
  \centering
  \begin{subfigure}[t]{0.33\textwidth}
    \centering
    \includegraphics[width=\linewidth]{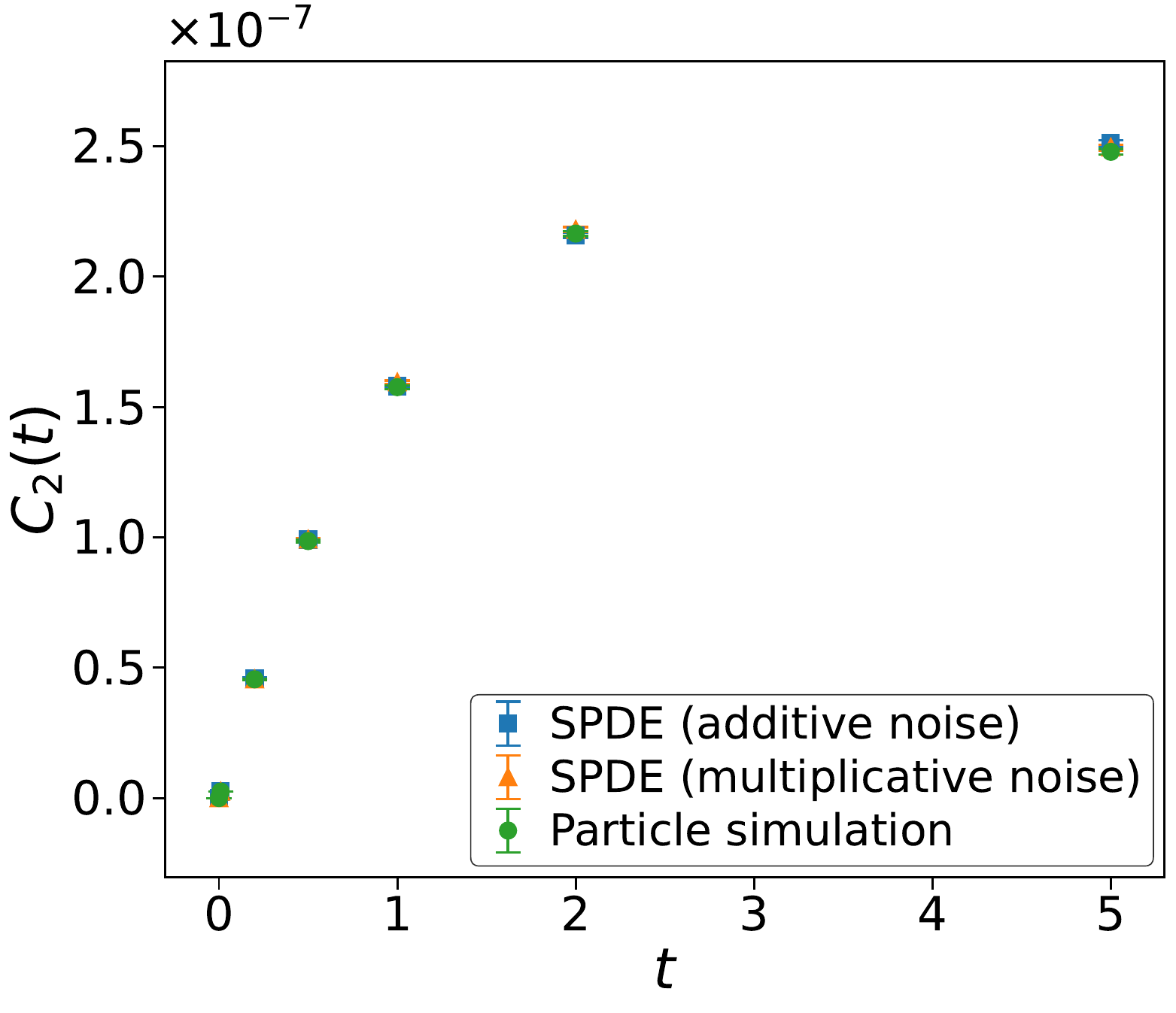}
    \subcaption{} 
    \label{fig:C2_sinx}
  \end{subfigure}\hfill
  \begin{subfigure}[t]{0.33\textwidth}
    \centering
    \includegraphics[width=\linewidth]{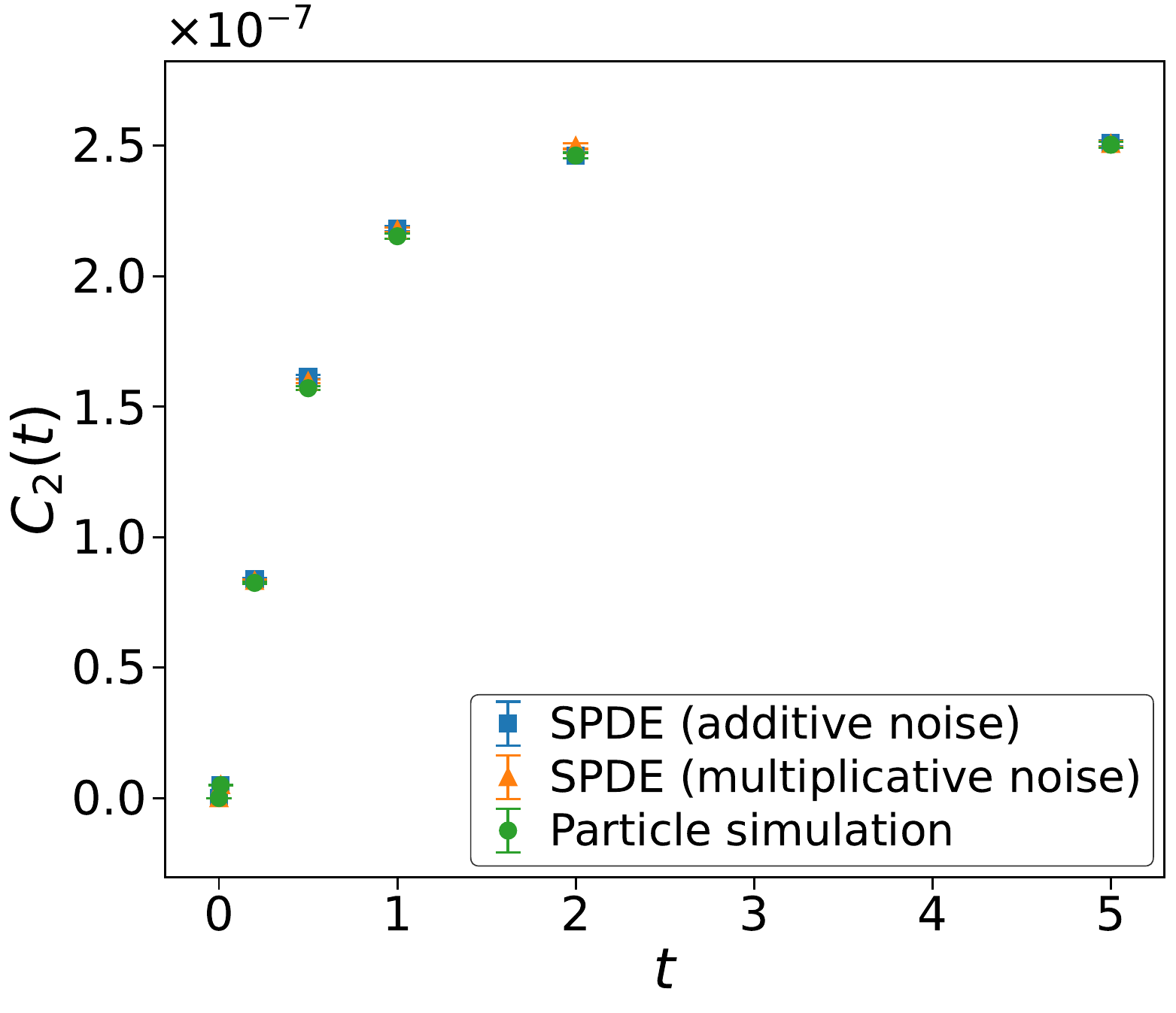}
    \subcaption{} 
    \label{fig:C2_sin2x}
  \end{subfigure}\hfill
  \begin{subfigure}[t]{0.33\textwidth}
    \centering
    \includegraphics[width=\linewidth]{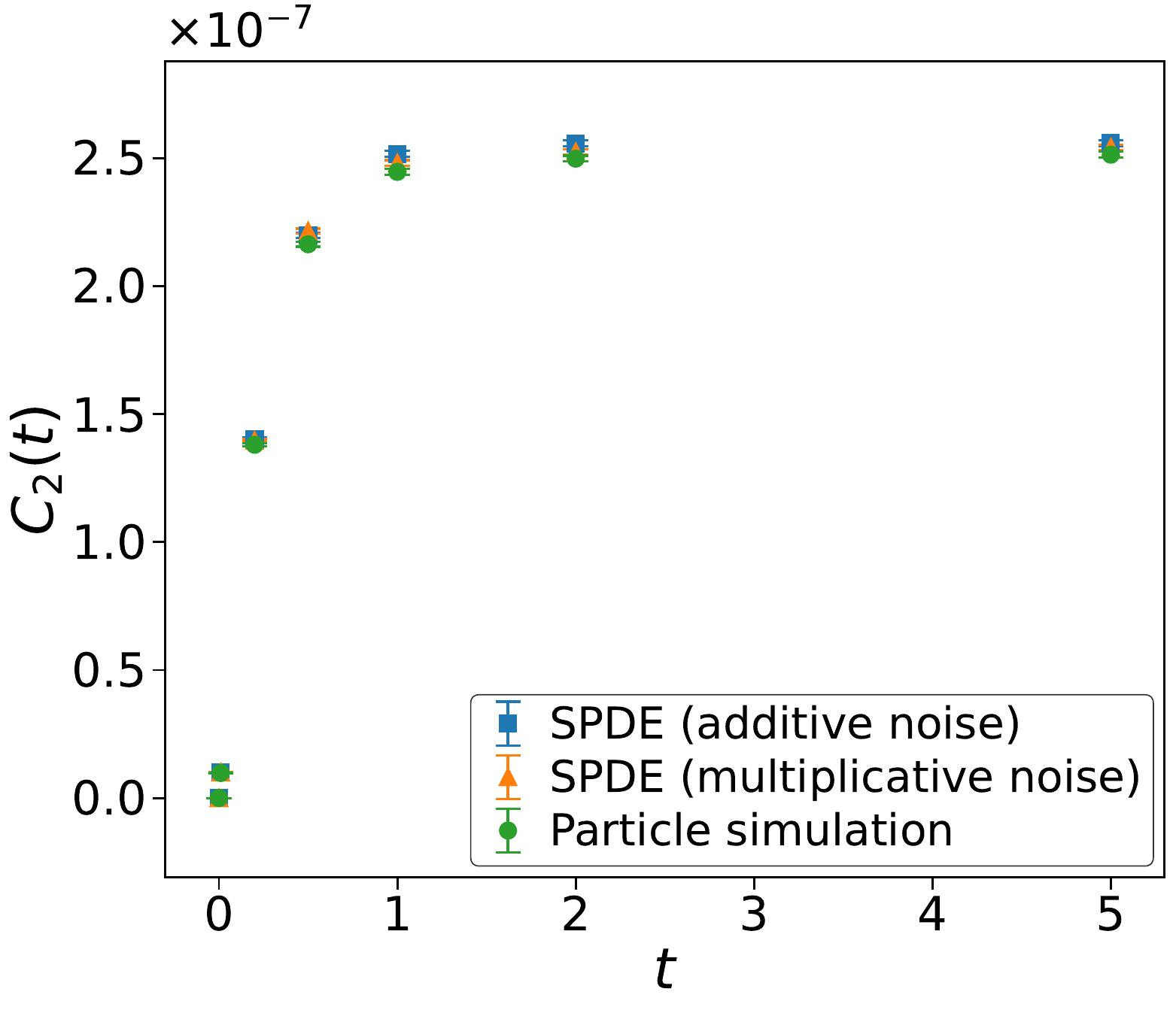}
    \subcaption{} 
    \label{fig:C2_sin4x}
  \end{subfigure}

  \caption{Time series of the variance $C_2(t)$ for $\varphi=\sin x$ (a), $\varphi=\sin 2x$ (b), and $\varphi=\sin 4x$ (c). The horizontal axis shows the physical time $t$ obtained from time steps with $dt=0.01$. Markers, with error bars, indicate the models: stars (blue) denote the SPDE with additive noise, triangles (orange) denote the SPDE with multiplicative noise, and circles (green) denote the particle simulation with $N=2\times10^6$, parameter $\alpha=1.0$, and $K(1)=1/2$. Symbols are almost superposed and error bars are barely visible.}
  \label{fig:C2_sin_panels}
\end{figure}
\captionsetup[subfigure]{labelformat=parens,labelsep=space}
\begin{figure}[tbp]
  \centering
  \begin{subfigure}[t]{0.48\textwidth}
    \centering
    \includegraphics[width=\linewidth]{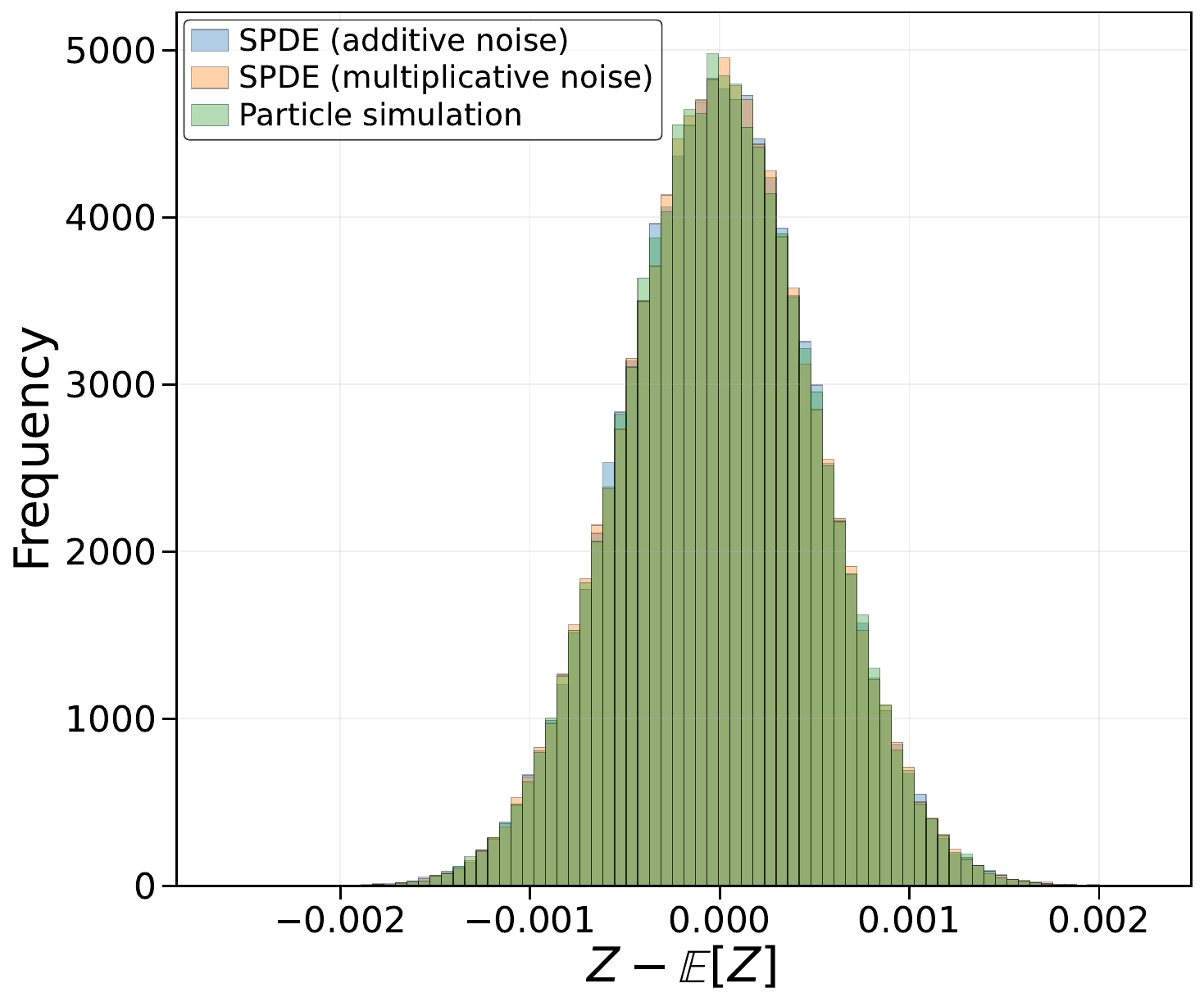}
    \subcaption{} 
    \label{fig:hist_sin2x_final}
  \end{subfigure}\hfill
  \begin{subfigure}[t]{0.48\textwidth}
    \centering
    \includegraphics[width=\linewidth]{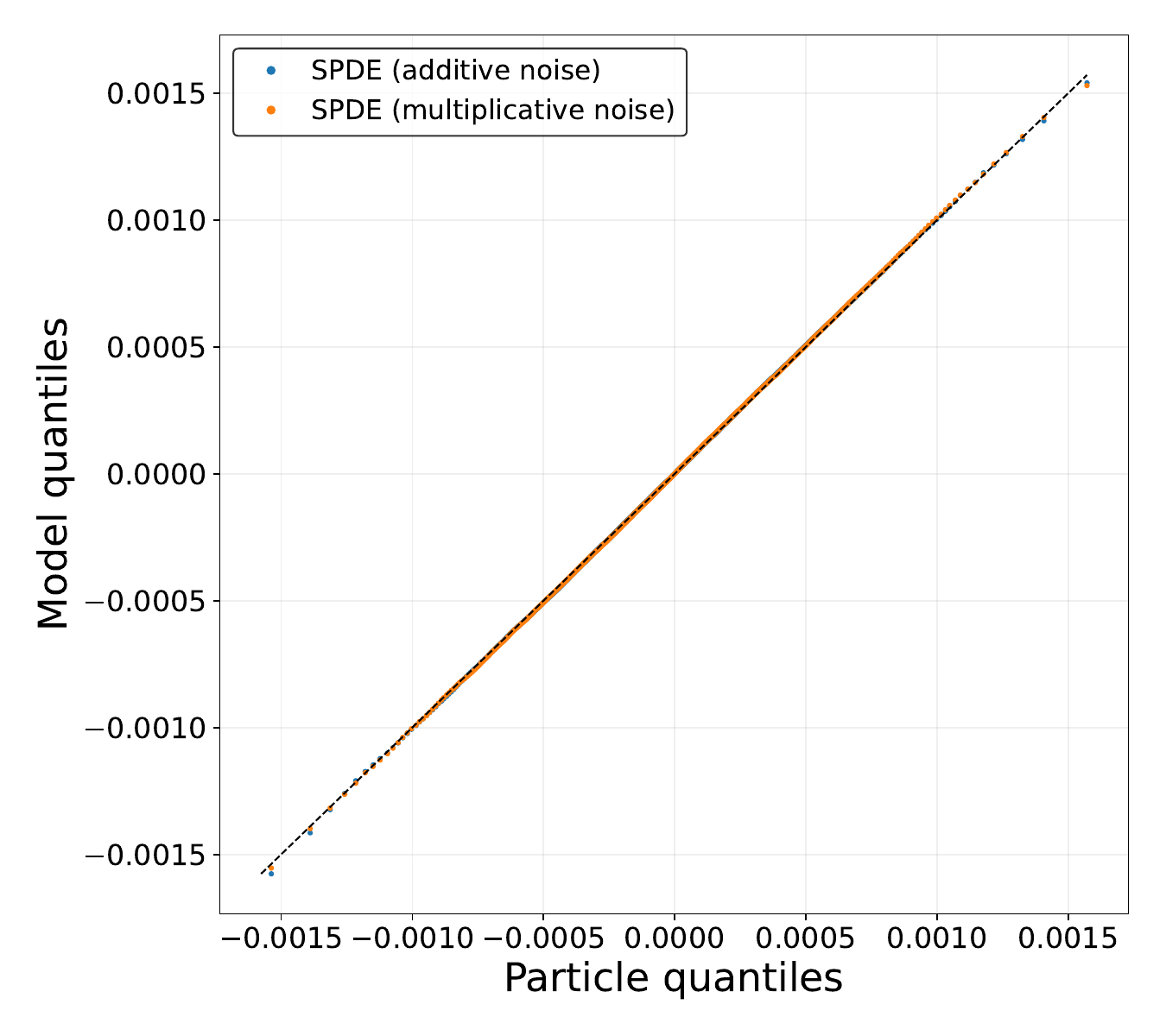}
    \subcaption{} 
    \label{fig:qq_sin2x_final}
  \end{subfigure}

  \caption{Distributional comparison for the representative test function $\varphi=\sin 2x$ at the final time step. The plotted quantity is the centered random variable $\langle \mu,\sin 2x\rangle-\mathbb{E}[\langle \mu,\sin 2x\rangle]$. Panel (a) shows the histogram of $\langle \mu,\sin 2x\rangle$, and panel (b) shows the QQ plot. Blue denotes the SPDE with additive noise, orange denotes the SPDE with multiplicative noise, and green denotes the particle simulation with $N=2\times10^6$, parameter $\alpha=1.0$, and $K(1)=1/2$.}
  \label{fig:dist_sin2x_final}
\end{figure}
In Fig.~\ref{fig:C2_sin_panels}, we plot the time series of $C_2(t)$ for the test functions $\varphi=\sin x$, $\varphi=\sin 2x$, and $\varphi=\sin 4x$, with the horizontal axis giving the physical time $t$. 
The same marker–color convention is used throughout: blue denotes the SPDE with additive noise, orange denotes the SPDE with multiplicative noise, and green denotes the particle simulation. 
Overall, both the additive-noise and multiplicative-noise models are in very good agreement with the particle-based simulation. Furthermore, once error bars are taken into account, the difference between the two SPDE models is very small. As shown in Fig.~\ref{fig:dist_sin2x_final}, the QQ plot suggests that the multiplicative-noise model may be slightly closer in the tails for the representative case $\varphi=\sin 2x$, but this effect remains weak. The quantitative accuracy we find indicates that the non-Gaussian effects encoded in the large deviation structure of section \ref{sec:ldp} are not visible in the probed regime.

\captionsetup[subfigure]{labelformat=parens,labelsep=space}
\begin{figure}[tbp]
  \centering
  \begin{subfigure}[t]{0.33\textwidth}
    \centering
    \includegraphics[width=\linewidth]{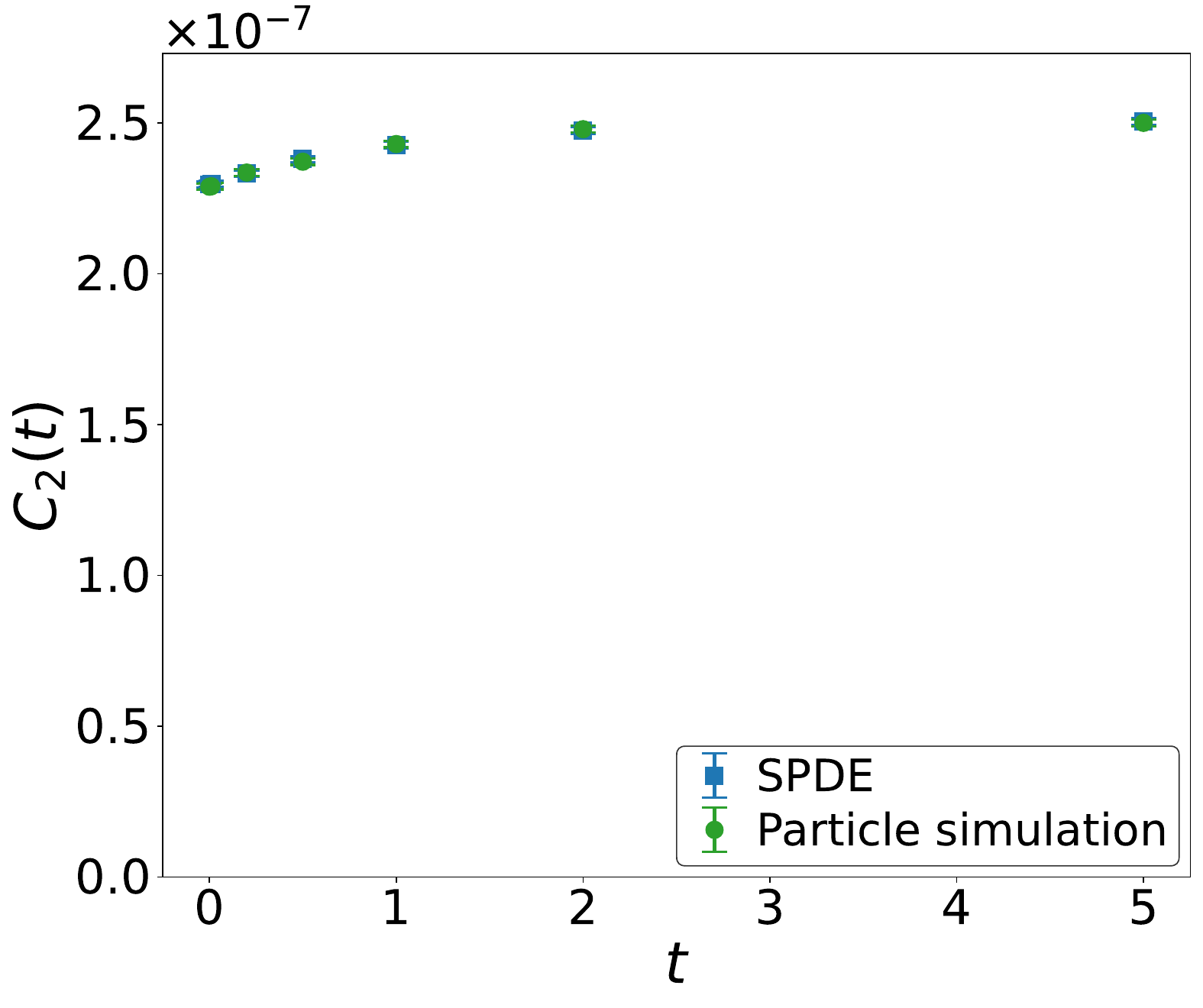}
    \subcaption{} 
    \label{fig:C2_sinx_random_initial}
  \end{subfigure}\hfill
  \begin{subfigure}[t]{0.33\textwidth}
    \centering
    \includegraphics[width=\linewidth]{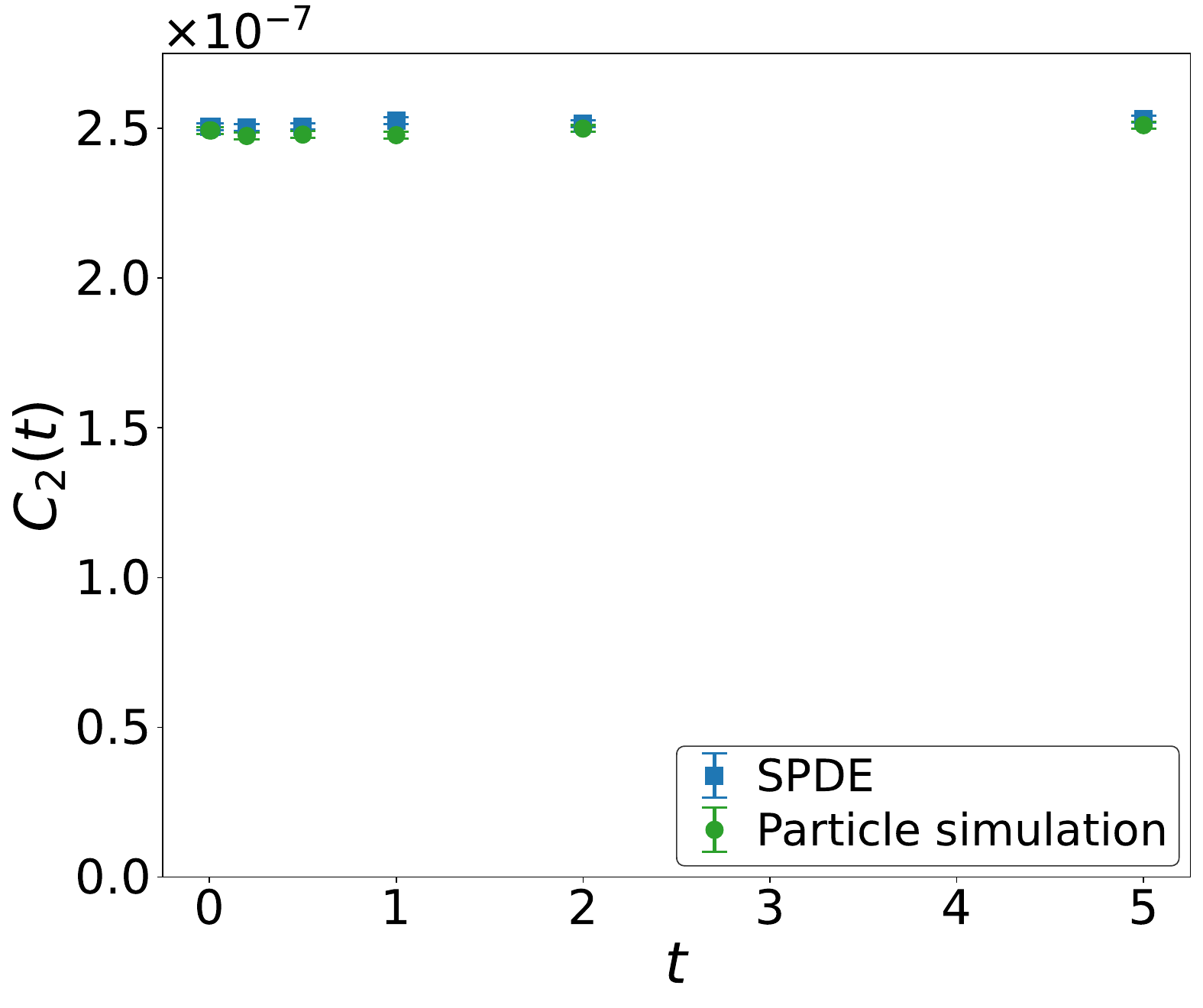}
    \subcaption{} 
    \label{fig:C2_sin2x_random_initial}
  \end{subfigure}\hfill
  \begin{subfigure}[t]{0.33\textwidth}
    \centering
    \includegraphics[width=\linewidth]{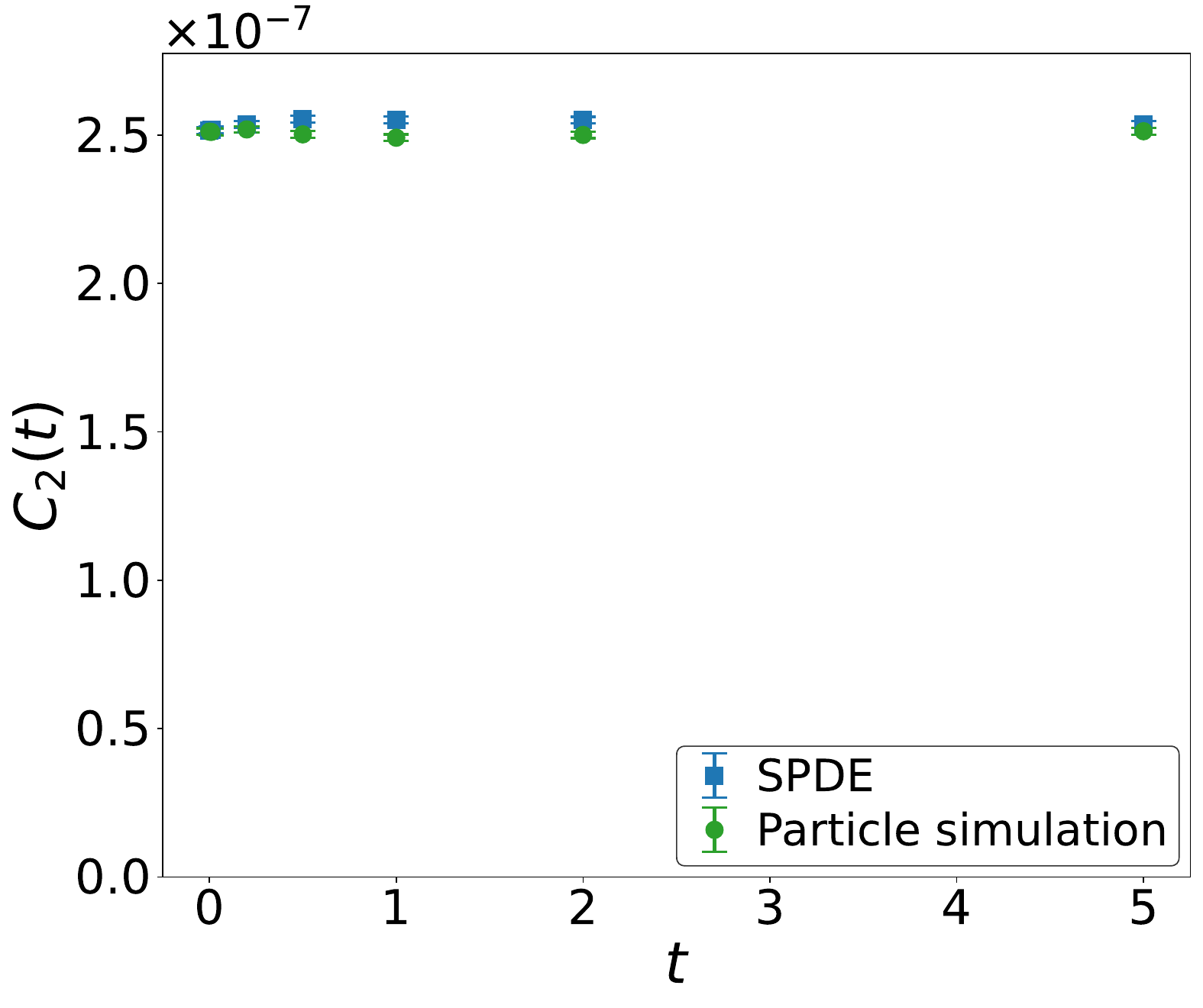}
    \subcaption{} 
    \label{fig:C2_sin4x_random_initial}
  \end{subfigure}

  \caption{Time series of the variance $C_2(t)$ for $\varphi=\sin x$ (a), $\varphi=\sin 2x$ (b), and $\varphi=\sin 4x$ (c) with random initial state. The horizontal axis shows the physical time $t$ obtained from time steps with $dt=0.01$. Markers, with error bars, indicate the simulations: stars (blue) denote the SPDE with additive noise, and circles (green) denote the particle simulation with $N=2\times10^6$, parameter $\alpha=1.0$, and $K(1)=1/2$.
  Symbols are almost superposed and error bars are barely visible.}
  \label{fig:C2_sin_panels_random_initial}
\end{figure}
\captionsetup[subfigure]{labelformat=parens,labelsep=space}
\begin{figure}[tbp]
  \centering
  \begin{subfigure}[t]{0.33\textwidth}
    \centering
    \includegraphics[width=\linewidth]{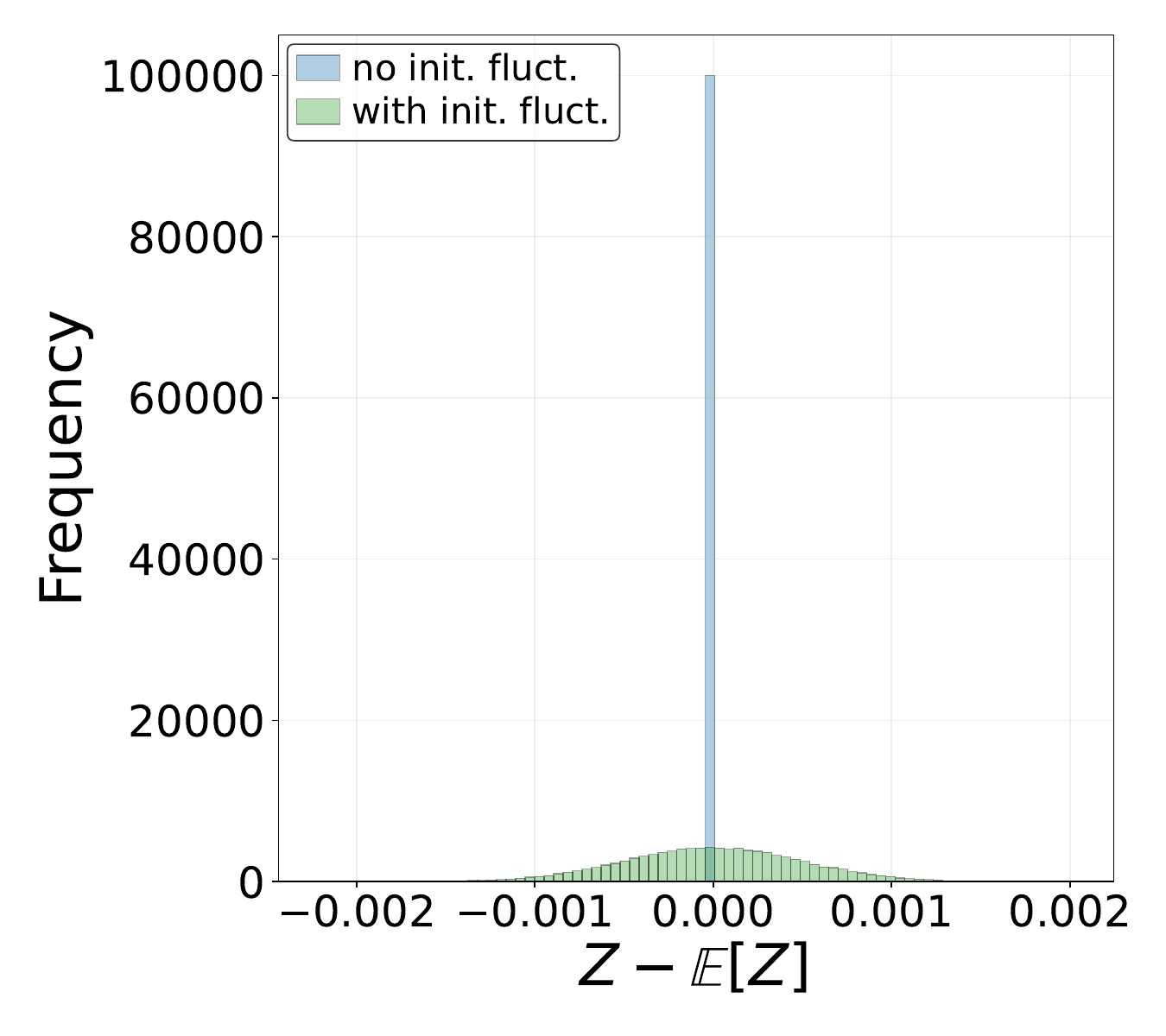}
    \subcaption{}
    \label{fig:hist_sin2x_self_initial_fluctuation_t0}
  \end{subfigure}\hfill
  \begin{subfigure}[t]{0.33\textwidth}
    \centering
    \includegraphics[width=\linewidth]{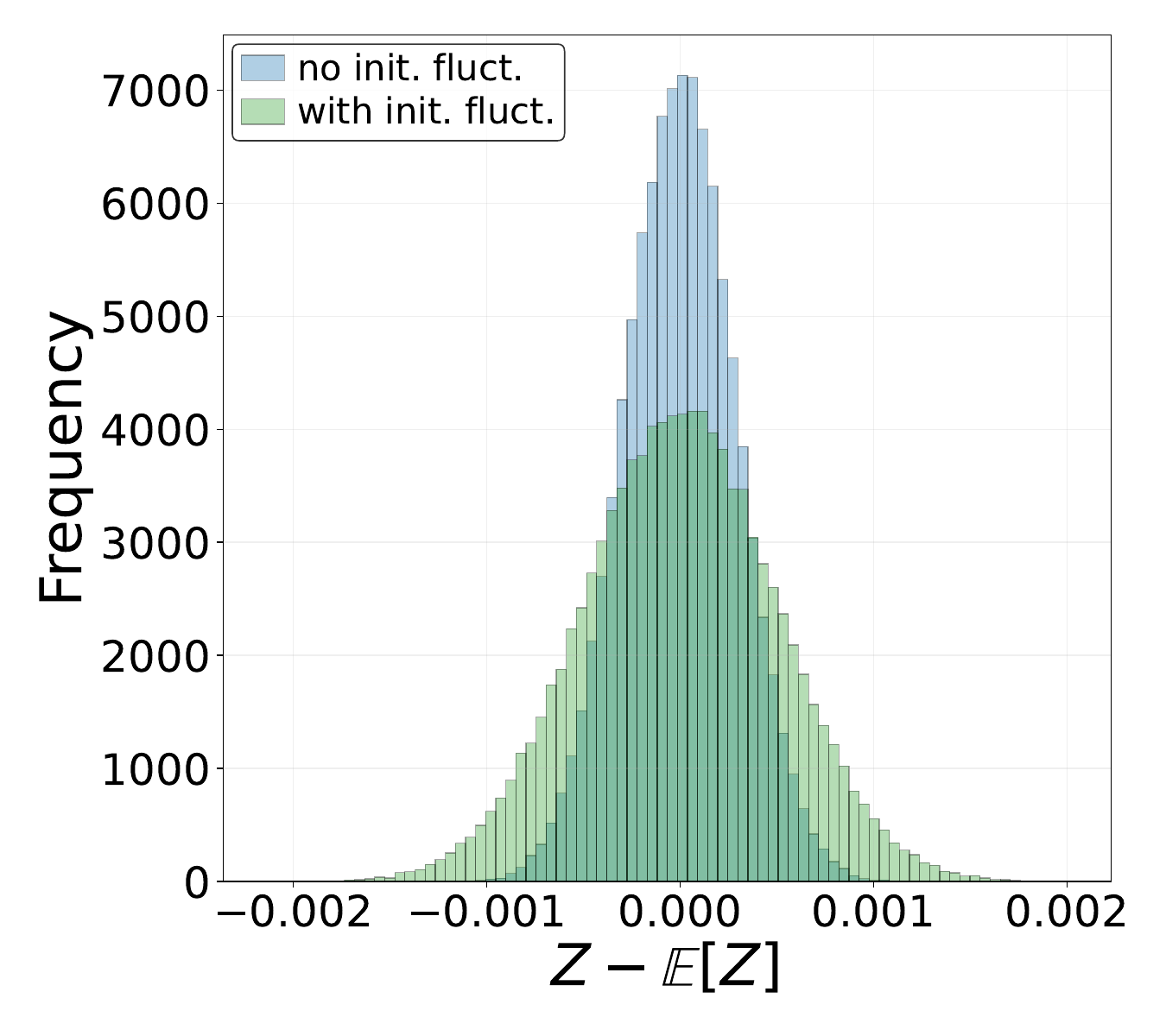}
    \subcaption{}
    \label{fig:hist_sin2x_self_initial_fluctuation_t02}
  \end{subfigure}\hfill
  \begin{subfigure}[t]{0.33\textwidth}
    \centering
    \includegraphics[width=\linewidth]{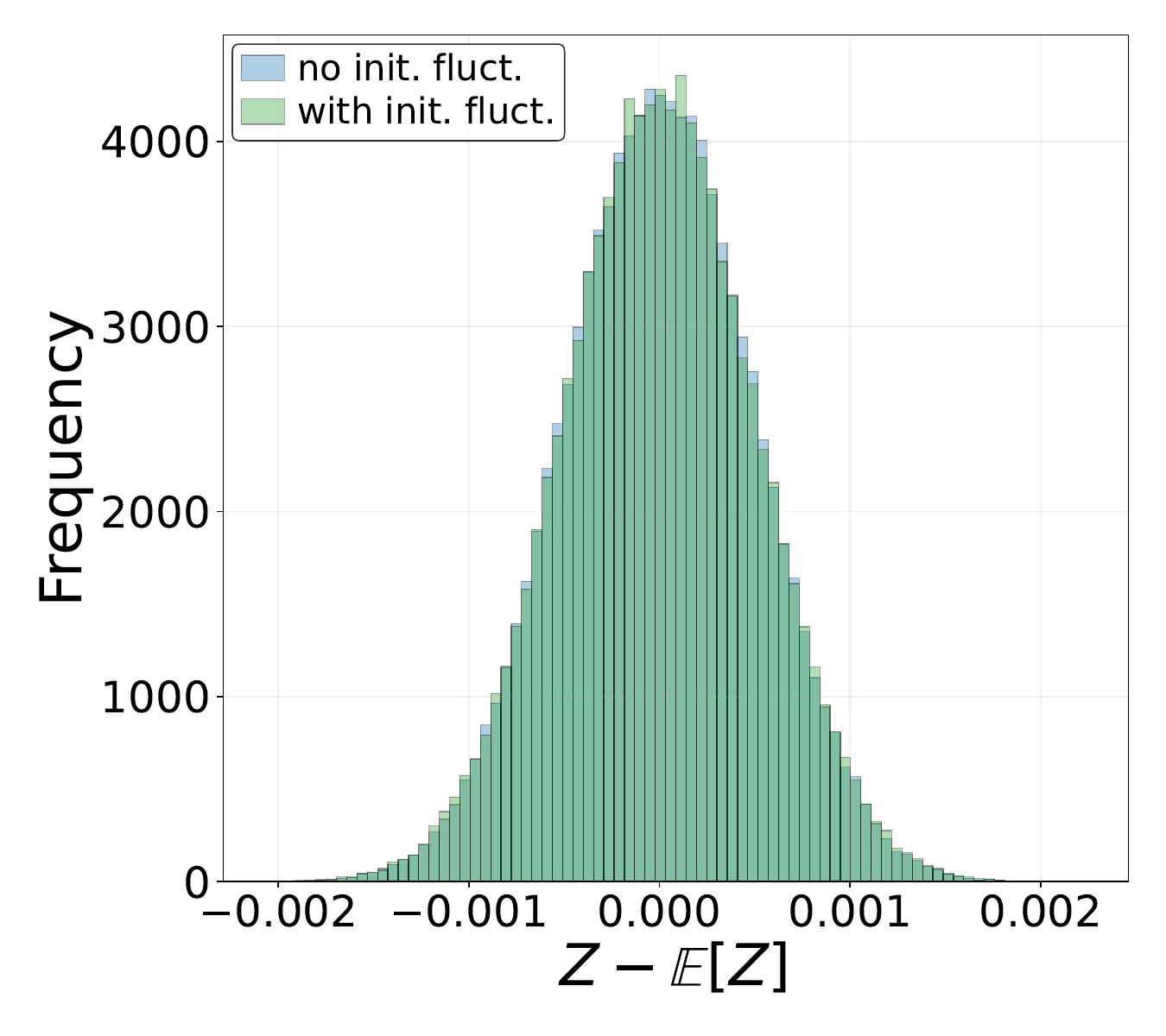}
    \subcaption{}
    \label{fig:hist_sin2x_self_initial_fluctuation_t5}
  \end{subfigure}

  \caption{Histogram comparison for SPDE model with and without initial fluctuation for the representative test function $\varphi=\sin 2x$. The plotted quantity is the centered random variable $Z-\mathbb{E}[Z]$, with $Z=\langle \mu,\sin 2x\rangle$. Panels (a), (b), and (c) correspond to $t=0$, $t=0.2$, and $t=5.0$, respectively. In each panel, blue denotes the case without initial fluctuation and green denotes the case with initial fluctuation.}
  \label{fig:hist_sin2x_self_initial_fluctuation}
\end{figure}

Instead of fixing the initial positions of the particles by the inverse CDF method, we now pick them up at random, independently, and with distribution \eqref{eq:IC}. 
In this case, fluctuations are already present at the initial time. Fig.~\ref{fig:C2_sin_panels_random_initial} shows the time series of the variance $C_2(t)$ for $\varphi=\sin x$, $\varphi=\sin 2x$, and $\varphi=\sin 4x$ when the initial state is random. As in the previous comparisons, the additive-noise SPDE \eqref{eq:TCL}, the multiplicative-noise SPDE \eqref{eq:DK_test}, and the particle simulation remain in relatively good agreement. Therefore, in the following, we show only the additive noise case in our plots. To generate the random initial condition for the SPDE, we proceed in a way consistent with the particle dynamics: for each realization, we construct an empirical histogram from the 
initial positions drawn at random for the microscopic simulation.
 This particle-based random profile is then used as the initial
condition for the SPDE. Consequently, the SPDE initial condition varies from one realization to another and contains the same finite-particle sampling fluctuations as the initial state of the particle simulation. With this initialization, the variance initially reflects the fluctuations already present at time $t=0$. However, after a sufficient time, the results approach those obtained without initial fluctuations, indicating that the system progressively loses the memory of the noisy initial state. This tendency is consistent with the histogram comparison in Fig.~\ref{fig:hist_sin2x_self_initial_fluctuation}, where the difference is clearest at $t=0$ and becomes much smaller at later times. 
\captionsetup[subfigure]{labelformat=parens,labelsep=space}
\begin{figure}[tbp]
  \centering
  \begin{subfigure}[t]{0.48\textwidth}
    \centering
    \includegraphics[width=\linewidth]{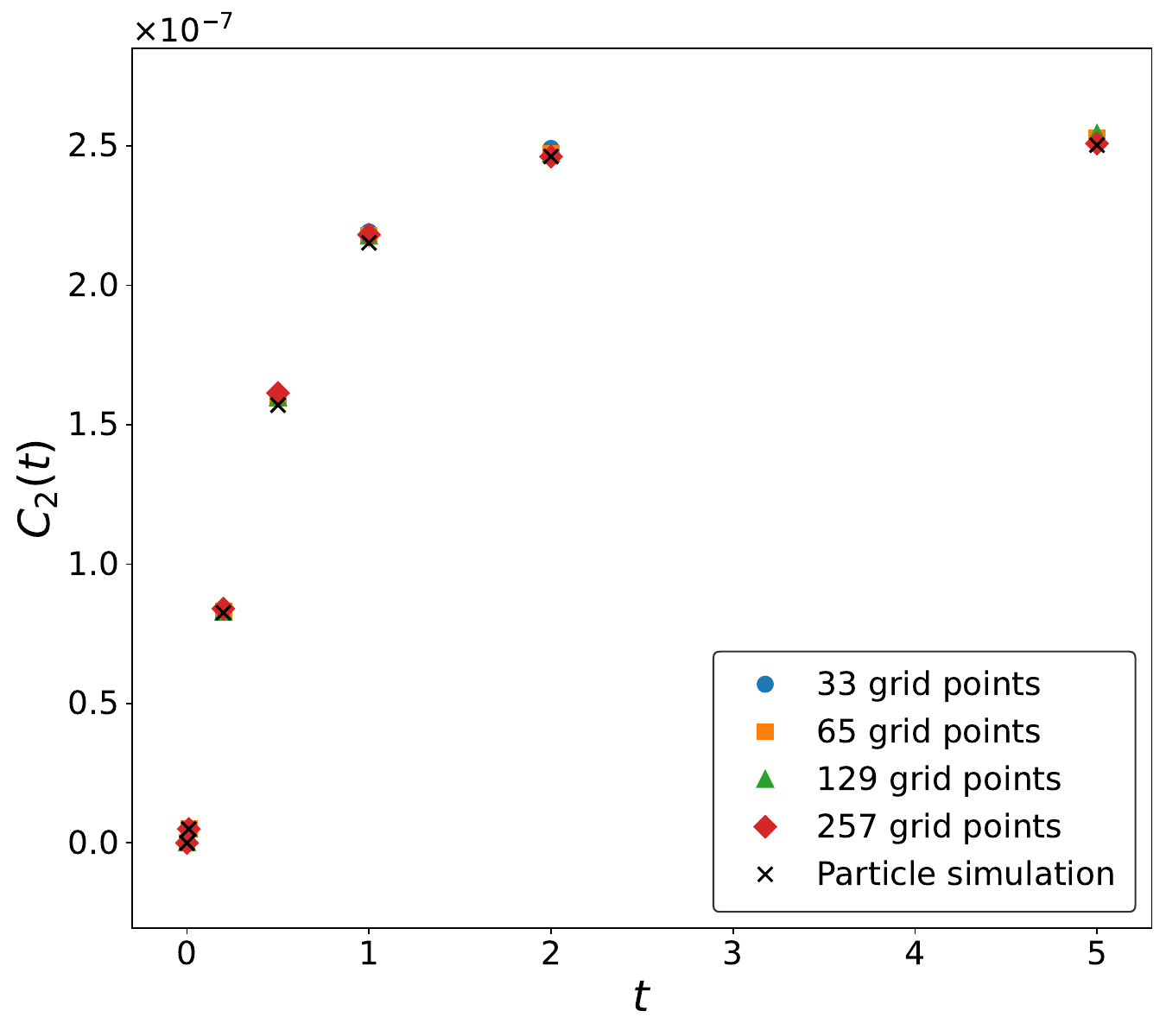}
    \subcaption{}
    \label{fig:hdep_c2_self_sin2x}
  \end{subfigure}\hfill
  \begin{subfigure}[t]{0.48\textwidth}
    \centering
    \includegraphics[width=\linewidth]{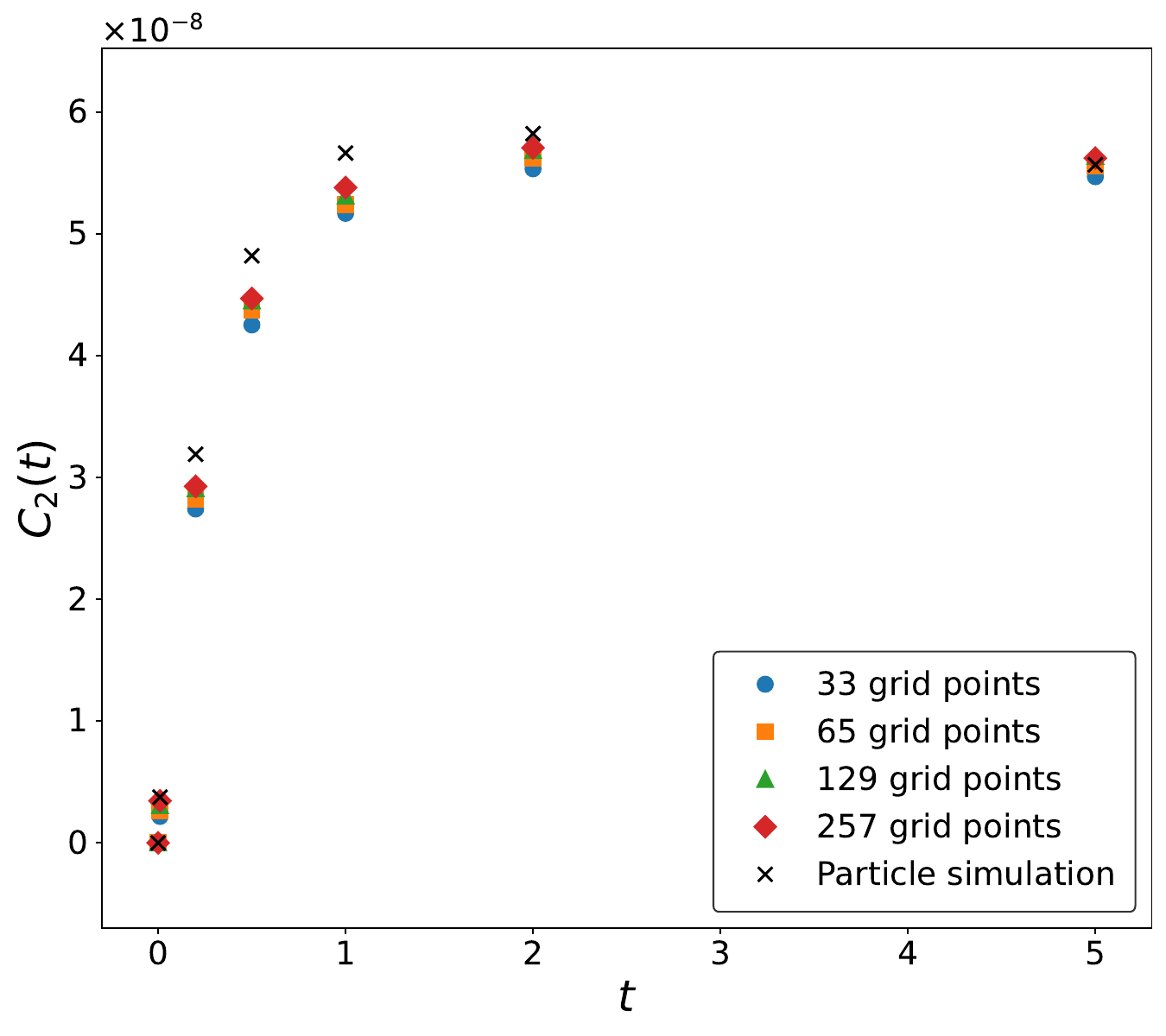}
    \subcaption{}
    \label{fig:hdep_c2_self_indicator}
  \end{subfigure}

\caption{
Grid size $h$ dependence of the second cumulant $C_2(t)$ for four spatial grids, $33$, $65$, $129$ and $257$, and comparison with the particle simulation. In both panels, star (blue) markers denote the result with $33$ grid points, squares (orange) denote $65$ grid points, triangles (green) denote $129$ grid points, diamonds (red) denote $257$ grid points, and crosses (black) denote the particle simulation. In all cases, the standard errors are approximately two orders of magnitude smaller than the corresponding values of $C_2$. Therefore, those are not shown for clarity. Panel (a) shows the result for the smooth test function $\varphi=\sin 2x$, where the four grid-based datasets almost overlap, and the dependence on the grid width is barely visible. Panel (b) shows the result for the non-smooth indicator test function supported on the interval $[L/2-L/16,\, L/2+L/16]$, where the differences among the grids become visible. This indicates that the spatial discretization effect is weak for smooth test functions, while it is clearer for non-smooth test functions.}
  \label{fig:h_dependence_c2}
\end{figure}
\begin{figure}[tbp]
  \centering
  \includegraphics[width=0.62\textwidth]{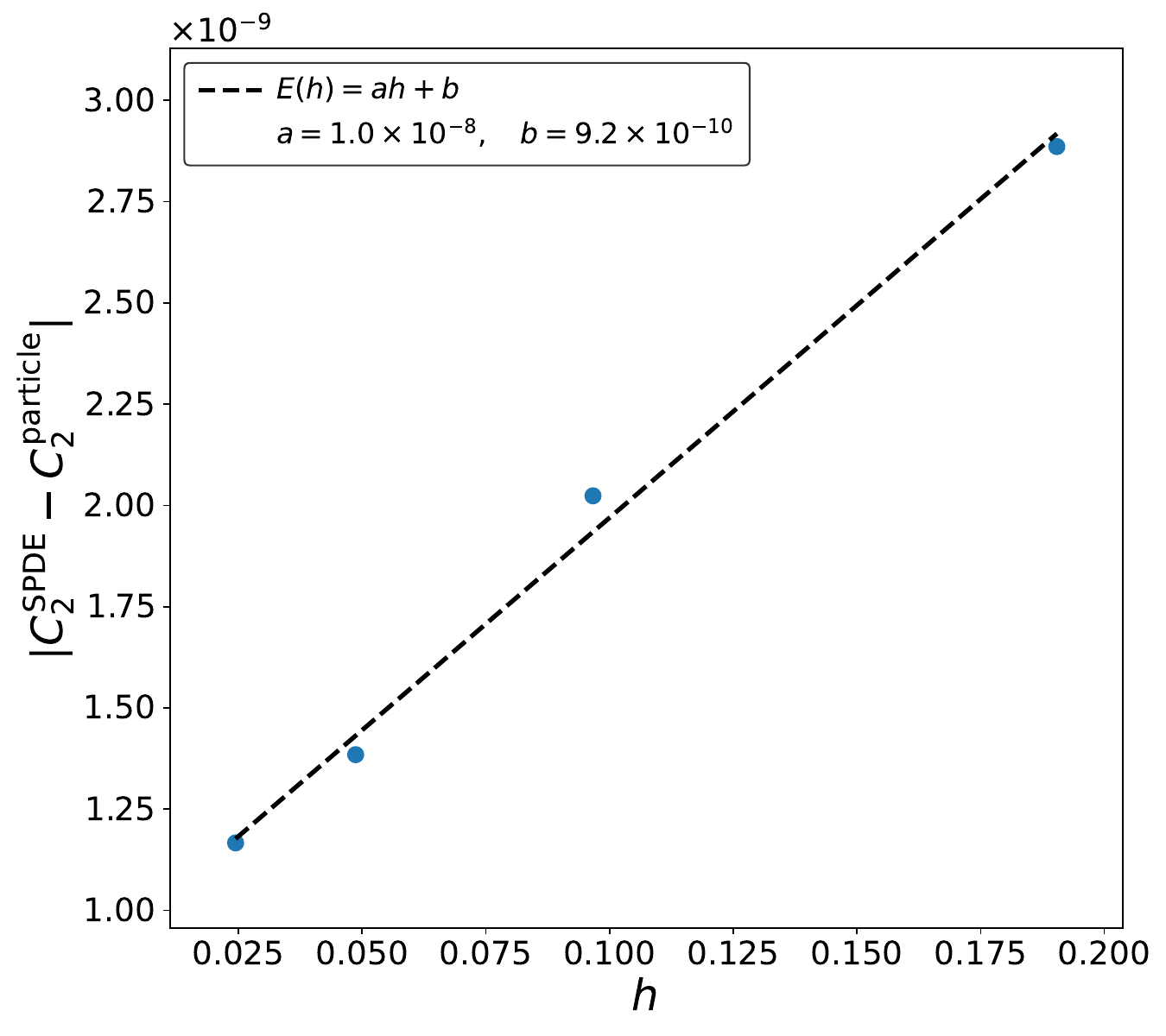}
  \caption{
  Dependence of the fixed-time discrepancy $E(h)=|C_2^{\mathrm{SPDE}}(h)-C_2^{\mathrm{particle}}|$
  on the spatial grid width $h$ at $t=2$ for the indicator observable supported on $[L/2-L/16,\,L/2+L/16]$. The additive-noise SPDE is compared with the particle simulation for four spatial grids, $33$, $65$, $129$, and $257$. The dashed line shows a linear fit of the form $E(h)=ah+b$. 
  }
  \label{fig:c2_particle_error_h_indicator}
\end{figure}

To evaluate the dependence on the spatial discretization, we compare the second cumulant $C_2(t)$ for several grid sizes. Fig.~\ref{fig:h_dependence_c2} shows the results obtained with $33$, $65$, $129$, and $257$ grid points. Panel (a) corresponds to the smooth test function $\varphi=\sin 2x$, while panel (b) corresponds to taking $\varphi$ as the indicator function of the interval $[L/2-L/16,\,L/2+L/16]$. For the smooth test function in Fig.~\ref{fig:hdep_c2_self_sin2x}, the dependence on the grid size is barely visible. By contrast, for the indicator function, the differences between the four grids become visible. In this sense, the effect of the spatial discretization remains weak for smooth observables, but becomes more apparent for a non-smooth observable. In addition, as shown in Fig.~\ref{fig:c2_particle_error_h_indicator}, we observe the discrepancy between the fractional SPDE and the particle simulation at the fixed time $t=2$ for the indicator observable. More precisely, we consider
\[
E(h)=\left|C_2^{\mathrm{SPDE}}(h)-C_2^{\mathrm{particle}}\right|,
\]
where $h$ denotes the spatial grid width of the SPDE discretization. Over the range of grid sizes considered here, this discrepancy decreases linearly with $h$. This indicates that the discrepancy observed for the indicator function is mainly due to the finite spatial resolution.

In conclusion, the numerical results indicate that the fractional SPDE \eqref{eq:TCL} captures the fluctuations of the particle system well at long times within the CLT regime, suggesting it provides an effective mesoscopic model for fluctuating fractional diffusion. The multiplicative noise variant \eqref{eq:DK_test} does not provide any significant improvement.
We also find that the effect of the initial fluctuations decreases as time evolves, and that the fluctuations at later times are mainly governed by the microscopic jump process itself. Finally, the dependence on the spatial discretization is weak for the smooth observables considered here, since they are mainly determined by low-frequency modes, which are accurately 
reproduced by our Fourier numerical method.
 However, it becomes more apparent for non-smooth observables, which are more sensitive to the truncation of high-frequency modes.

\section{Conclusion}
To summarize, the main results of this work are:
\begin{itemize}
\item The derivation of a Large Deviation Principle for the empirical measure of particles submitted to Lévy noise. In the case of independent particles, this unveils a generalized gradient flow structure.
\item The derivation of the linear fluctuating hydrodynamics equation in this setting, with explicit expressions for the long-range correlations of the noise.
\item A numerical implementation of the linear fluctuating hydrodynamics, showing very good agreement with direct particle simulations, always in the independent case.  
\end{itemize}
These results, mostly obtained in a simplified setting, 
open the door to many further studies. At the theoretical level, one may wonder whether the $\alpha-$stable McKean-Vlasov processes of section \ref{sec:mean-field-levy} have a gradient flow structure, and what it would be. It is also noteworthy that the structure of the LDP for the empirical measure of Lévy particles is different from the one found in \cite{Bernardin2025} starting from a long-range interacting Ginzburg-Landau model, whose hydrodynamic limit is a fractional diffusion. This raises the question for other microscopic models behaving at large scales as a fractional diffusion: are their fluctuations locally Gaussian as in \cite{Bernardin2025} and \cite{bernardin_etal_2026}, or "jump-like" as in the present paper? One could think, for instance, of studying the question for microscopic models whose kinetic formulation has been shown to converge to a fractional diffusion in the hydrodynamic limit \cite{jara2009limit,Mellet2010,Mellet2011,Abdallah2011,Abdallah2011a,cattiaux2014limit}.
At the numerical level, the implementation of the linear fluctuating hydrodynamics proposed in this paper should extend straightforwardly to the mean-field interacting case. Going beyond the linear setting, to include the full LDP structure is much more challenging.

\section*{Acknowledgments}
The authors gratefully acknowledge discussions with C. Bernardin, R. Chetrite and M. Puel. This work is partially supported by the grants PERISTOCH ANR-19-CE40-0023 and RETENU
ANR-20-CE40-0005-01 of the French Agence Nationale
de la Recherche. H.H. also acknowledges the support of the Agence Nationale de la Recherche through the MilaDy project, Grant No. ANR-25-CE51-0077-02.

\clearpage
\appendix           
\section{On the fractional Laplacian}
\label{app:frac_lapl}

We gather here three equivalent definitions of the fractional Laplacian, with references.

{\bf Spectral definition:} It may be the simplest, cf. the main text. 

With $\Omega=[0,L]$ and periodic boundary conditions, it reads as follows. A function $f$ defined on $\Omega$ with periodic boundary conditions can be seen equivalently as an $L$-periodic function defined on $\mathbb{R}$. Then for any real $x$:
\[
f(x) =\sum_{k\in \mathbb{Z}}c_k(f) e^{2i\pi k x/L},
\]
where the $c_k(f)$ are the Fourier coefficients of $f$.
Then
\[
-(-\Delta)^{\alpha/2}f (x) =\sum_{k\in \mathbb{Z}} -\left|\frac{2\pi k}{L}\right|^\alpha c_k(f)e^{2i\pi k x /L}.
\]

{\bf Generator definition:}\\
Let $X_\alpha(t)$ be an $\alpha-$stable symmetrical Levy process in $\mathbb{R}$ (it is unique up to normalization), normalized such that for any $t$
\[
\mathbb{E}\left[e^{ikX_\alpha(t)}\right] =e^{-t|k|^\alpha}.
\]
Then the generator of the process $X_\alpha$ is exactly the fractional Laplacian: for any test function $u$
\[
\lim\limits_{h\to 0} \frac{\mathbb{E}[u(x+X_\alpha(h))-u(x)]}{h} = -(-\Delta^{\alpha/2})u(x).
\]
Consider now the periodized process $\tilde{X}_{\alpha}(t)=X_\alpha(t) [L]$ on $[0,L]$ with periodic boundary conditions. Let $u$ be a test function ($u$ is $L-$periodic). Then 
\[
u(x+\tilde{X}_\alpha(h))=u(x+X_\alpha(h)),
\]
and a straightforward computation shows that
\[
\lim\limits_{h\to 0} \frac{\mathbb{E}[u(x+\tilde{X}_\alpha(h))-u(x)]}{h} = -(-\Delta^{\alpha/2})u(x),
\]
where the spectral definition of the fractional Laplacian on $[0,L]$ with periodic boundary conditions is used. 

{\bf Singular integral definition:}\\
From \cite{lischke2020fractional}: 
\begin{equation}
(-\Delta)^{\alpha/2}f(x) = C(\alpha)\, {\rm p.v.}\!\int_{\mathbb{R}} \frac{f(x)-f(y)}{|x-y|^{1+\alpha}}dy,
\end{equation}
with
\[
C(\alpha) = \frac{2^\alpha \Gamma\left(\frac{1+\alpha}{2}\right)}{\sqrt{\pi}|\Gamma\left(-\frac{\alpha}{2}\right)|}.
\]
\cite{Mellet2010} does not give the expression of the constant $C(\alpha)$.

\section{Computation of the noise covariance matrix}
\label{app:calc}
Let us consider  $\boldsymbol{\Sigma}_{\mathbf{xx}}$ as a computation example. First of all, we should calculate all correlation functions with complex conjugates of the Fourier coefficients of the noise. Recalling the noise correlation \eqref{eq:noise_fourier}, those are given by
\[
  \begin{aligned}
    \mathbb{E}[\hat{\xi}(k_1)\hat{\xi}(k_2)]
      &= C(k_1,k_2),\\
    \mathbb{E}[\hat{\xi}(k_1)\hat{\xi}^*(k_2)]
      &=  C(k_1,-k_2),\\
    \mathbb{E}[\hat{\xi}^*(k_1)\hat{\xi}(k_2)]
      &=  C(-k_1,k_2),\\
    \mathbb{E}[\hat{\xi}^*(k_1)\hat{\xi}^*(k_2)]
      &=  C(-k_1,-k_2).
  \end{aligned}
\]
Meanwhile, $\hat{\xi}$ is $\hat{\xi}_{\rm Re}+i\hat{\xi}_{\rm Im}$, then we can obtain the decomposition formula such as 
\[
  \begin{aligned}
    \mathbb{E}[\hat{\xi}(k_1)\hat{\xi}(k_2)]
    &= \mathbb{E}[\hat{\xi}_{\rm Re}(k_1)\hat{\xi}_{\rm Re}(k_2)]
      - \mathbb{E}[\hat{\xi}_{\rm Im}(k_1)\hat{\xi}_{\rm Im}(k_2)] \\
    &\quad{}+\,i\,\mathbb{E}[\hat{\xi}_{\rm Re}(k_1)\hat{\xi}_{\rm Im}(k_2)]
      +\,i\,\mathbb{E}[\hat{\xi}_{\rm Im}(k_1)\hat{\xi}_{\rm Re}(k_2)].
  \end{aligned}
\]
Using the above decompositions, the block components of the covariance matrix $\Sigma$ can be obtained.  For example, the $\Sigma_{xx}$ block is given by:
\[
  \begin{aligned}
    \boldsymbol{\Sigma}_{\mathbf{xx}}
    = \mathbb{E}[\hat{\xi}_{\rm Re}(k_1)\hat{\xi}_{\rm Re}(k_2)] 
    = \frac{1}{4}\bigl(
       \mathbb{E}[\hat{\xi}(k_1)\hat{\xi}(k_2)]
      +\mathbb{E}[\hat{\xi}(k_1)\hat{\xi}^*(k_2)] \\
    \quad\quad{}
      +\mathbb{E}[\hat{\xi}^*(k_1)\hat{\xi}(k_2)]
      +\mathbb{E}[\hat{\xi}^*(k_1)\hat{\xi}^*(k_2)]
    \bigr).
  \end{aligned}
\]
Therefore, we can obtain $\boldsymbol{\Sigma}_{2n \times 2n}$ to generate this complex random variable.

\bibliographystyle{unsrtsiamplain}
\bibliography{references}
\end{document}